\documentclass[journal,final]{IEEEtran}
\usepackage[utf8]{inputenc}
\usepackage{graphicx, flushend}
\usepackage{amsmath}
\usepackage{cases}
\usepackage{bbm}
\usepackage{amssymb}
\usepackage{array}
\usepackage{lipsum,color}
\usepackage{epstopdf}

\DeclareGraphicsExtensions{.pdf,.png,.jpg,.eps,.ps}

\usepackage{algorithm}
\usepackage{algpseudocode}
\usepackage[square, numbers, comma, sort&compress]{natbib}
\usepackage{xcolor}
\usepackage{hyperref}

\newcommand{\bs}{\boldsymbol}

\usepackage{url}
\usepackage{amsthm,amssymb}

\newtheorem{definition}{Definition}

\newcommand{\ba}{\begin{array}}
	\newcommand{\ea}{\end{array}}

\usepackage{textcomp}
\usepackage[font=small,skip=4pt]{caption}

\begin{document}
	

\title{High-Precision Intelligent Reflecting Surfaces-assisted Positioning Service in 5G Networks with Flexible Numerology}
\pagestyle{empty}

\author{\IEEEauthorblockN{Ti Ti Nguyen, Kim-Khoa Nguyen}  \\ {\em Ecole de Technologie Sup\'{e}rieure,   Canada} }
	
 


\maketitle
\thispagestyle{empty}
\begin{abstract}
Accurate positioning is paramount for a wide array of location-based services (LBS) in fifth-generation (5G) wireless networks. Recent advances in 5G New Radio (NR) technology holds promise for very high-precision positioning services. Yet, challenges arise due to diverse types of numerology and massive connected devices. This paper presents a novel approach to improve positioning precision within a 5G NR framework with comb patterns on time-frequency resource mapping. We then formulate an optimization problem  aimed at minimizing the maximum users' positioning error in an intelligent reflected surface (IRS)-assisted 5G network by controlling the user-anchor association, numerology-related selection, IRS's reflecting elements, privacy protection level, and transmit power. To address the non-convex nature of the underlying mixed-integer non-convex problem (MINLP), we propose an efficient algorithm that combines optimization, matching, and learning techniques. Through extensive numerical experiments, we demonstrate the effectiveness of our proposed algorithm in minimizing positioning errors compared to conventional methods.

\end{abstract}
%
\section{Introduction}

Nowadays, the location-based services (LBS) are emerging as integral components in diverse sectors like logistics, smart factories, autonomous vehicles, localized sensing, asset tracking, digital twins, augmented and virtual reality \cite{guo2019survey, del2017survey}. Traditional LBS relying on global navigation satellite systems (GNSS) are often limited to urban environments  due to signal obstruction and disruption caused by buildings  \cite{li2022carrier}. 
Therefore, communication-based positioning methods, such as 5G, Bluetooth, and Wi-Fi, have experienced significant progress in recent years. Leveraging existing communication networks, these methods offer inherent advantages, including extensive coverage, easy deployment, and cost-effectiveness \cite{yin2022design}. 


 Within communication-based positioning methods, 
the work in \cite{wei20225g} utilizes the positioning reference signal (PRS) of 5G to evaluate the positioning error of a single-user system. The study demonstrates that the range error decreases over 7.5 times when applying the PRS with 256 subcarriers compared to the synchronization signal (SS) with 127 subcarriers. Additionally, the range error increases with larger comb sizes, although this difference gradually diminishes with the improvement of signal-to-noise ratio (SNR).
In \cite{pan2022efficient}, the authors investigate the joint angle and delay estimation (JADE) problem in a multiple-input multiple-output (MIMO) network and propose an algorithm with lower complexity than the multiple signal classification (MUSIC)-based JADE method \cite{kotaru2015spotfi}. Moreover, in \cite{driusso2016vehicular}, the authors apply the ESPRIT and Kalman filter for time-of-arrival tracking (EKAT) algorithm to tracking the time-of-arrival (ToA) of a vehicle. 
For multi-user systems, traditional methods primarily focused on orthogonal resource allocation for positioning services among users. However, the surge in connected mobile devices within upcoming B5G and 6G networks presents significant challenges to this conventional approach  \cite{yu2022location, yin2022design}.

The non-orthogonal resource allocation approach has emerged as a solution to address limited wireless resources. However, in 5G NR, where flexible numerology allows the application of different sub-carrier spacings (SCS) \cite{nguyen2019wireless, ha2019admission}, managing interference becomes crucial for efficient utilization of wireless resources \cite{martin2021emerging}. In \cite{ghimire2022reference}, different cyclic shifts to Zadoff-Chu sequences (ZC) are applied to reduce multi-user interference in time difference of arrival (TDoA) estimation.
Additionally, the impact of interference on ToA estimation error is analyzed in \cite{betz2009generalized} based on code tracking with an early discriminator design. 
Recently, advances in the intelligent reflecting surface (IRS) technology have promised advantages in manipulating the spatial distribution of wireless signals \cite{tang2022physical}. 
Beyond communications, IRS applications have extended to wireless localization, as explored in previous studies \cite{zhang2023multi, yu2022location, hu2018beyond, zhang2021metalocalization}.
In \cite{zhang2023multi}, a framework is investigated to integrate IRS technology, breaking down the spatial resolution limitation of WiFi signals in both single-person and multi-person scenarios. In \cite{yu2022location}, a multi-user location sensing algorithm is developed based on uplink communication signals transmitted from multiple users to the base station (BS). 
In \cite{hu2018beyond}, Cramer-Rao Lower Bounds (CRLB) are derived for positioning with IRS. 
In \cite{zhang2021metalocalization}, IRS is applied to create a favorable received signal strength (RSS) distribution. The study shows that the localization error is reduced by at least threefold.

In communication-based positioning methods, with the increasing amount of data collected and transmitted over these networks, privacy concerns have become a major issue for individuals and organizations \cite{nguyen2021security}. Differential privacy (DP) has emerged as a promising solution for safeguarding data privacy in wireless communications \cite{liu2020privacy}. By ensuring that statistical analysis of data does not reveal sensitive information about individuals, differential privacy mathematically guarantees that an individual's private information cannot be inferred during data processing. This is achieved by adding noise to the data to ensure the statistical similarity of the output analysis, regardless of whether a particular individual's data is included or not \cite{geng2015optimal, ligett2017accuracy}.
Recent studies have focused on the application of differential privacy in positioning systems \cite{huang2019incentivizing, zhao2020local, kim2019workload}.  In \cite{huang2019incentivizing}, a Bayesian game is adapted to implement differential privacy in location-based crowd sensing services, which can restrict the sensitive location information leakage. 
The scheme proposed in \cite{zhao2020local} utilizes gradient perturbation techniques to ensure local differential privacy of data transmitted in internet of vehicle (IoV) systems. The work in \cite{kim2019workload}  explores a data encoding and perturbation scheme of local differential privacy (LDP) to minimize the overall error in estimating crowd densities in indoor spaces. Furthermore, in \cite{wei2019differential}, a differential privacy-based location protection (DPLP) scheme is introduced to protect the location privacy of both workers and tasks by splitting the exact locations of both workers and tasks into noisy multi-level grids.

However, none of aforementioned works  \cite{martin2021emerging, wei20225g, pan2022efficient,kotaru2015spotfi,driusso2016vehicular, ghimire2022reference, yin2022design, betz2009generalized, zhang2023multi, yu2022location,hu2018beyond,zhang2021metalocalization} have ever addressed the challenges of positioning in upcoming B5G and 6G networks with strict requirements on the privacy and interference management. In this paper,  we fill this research gap by modeling the interference due to non-orthogonal resource allocation  in the context of wireless networks with flexible numerology and analyzing the ToA estimation error in 5G NR systems with a comb pattern of positioning reference signals.

Our contribution is summarized as follows:
\begin{itemize}
	\item We present  a novel model for positioning services in 5G NR, taking into account transmitted orthogonal frequency-division multiplexing (OFDM) symbols with comb patterns and various numerology types. Then, we conduct a comprehensive analysis of interference due to the non-orthogonal resource allocation and derive the ToA error between anchors and users.
	\item We employ the differential privacy criteria to safeguard the location information of anchors while they provide communication-based positioning services to users.
	\item We formulate an optimization model to minimize the maximum users' positioning error in an IRS-assisted 5G network by controlling the user-anchor association, numerology-related selection, IRS's reflecting elements, and transmit power.
	\item We develop an effective algorithm that combines optimization, matching, and learning techniques to address the underlying mixed-integer non-linear programming (MINLP) problem. This approach optimizes continuous variables such as transmit power and privacy protection levels. Numerology and offset selections in a comb pattern are determined through many-to-one matching. User-anchor associations are achieved through many-to-many matching. Additionally, deep Q-network (DQN) learning is utilized for beam selection of IRS reflecting elements. 
	\item  We conduct extensive experiments to evaluate the performance of our proposed designs. The results validate the superior performance of our proposed algorithm.
\end{itemize}
The remainder of this paper is organized as follows.
Section \ref{sec_2} presents the system model, the positioning error analysis in 5G NR with the comb pattern of positioning reference signals. Section \ref{sec_4} introduces the proposed algorithm to address the underlying MINLP problem. Section~\ref{sec_5} evaluates the performance of the proposed algorithm. Finally, Section \ref{sec_6} concludes the paper and presents future directions.

\section{System model}
\label{sec_2}
Considering an indoor positioning system (IPS) including an IRS, $J$ anchors, and $K$ users  (referred to as `users' for simplification).  Let $\mathcal{J}$ and $\mathcal{K}$ be the set of anchors and users, respectively. Among the $J$ anchors, one is designated as an access point (AP) and possesses precise location information.
It is assumed that the remaining anchors are equipped with highly accurate location estimation sensors, such as LiDAR, D-RGB, and UWB. The location of anchor $j$ is determined by these sensors and can be expressed as follows
\begin{equation}
	\bs{s}_{j, \sf a} = \bs{\hat{s}}_{j, \sf a} + \bs{\eta}_{j, \sf a},
\end{equation}
where $\bs{\eta}_{j, \sf a} \sim \mathcal{N}(\bs{0}, \xi_{j,0}^2 \bs{I}_3)$, $\bs{s}_{j, \sf a}$ and $\bs{\hat{s}}_{j, \sf a}$ are the true and estimated location of the anchor $j$. Note that $\bs{s}_{j, \sf a} = [s_{j, \sf a,x}, s_{j, \sf a,y}, s_{j, \sf a,z}]$ and $\bs{\hat{s}}_{j, \sf a} = [\hat{s}_{j, \sf a,x}, \hat{s}_{j, \sf a,y}, \hat{s}_{j, \sf a,z}]$. We also denote the location of  device $k$ as $\bs{s}_{k, \sf r} = [s_{k, \sf r,x}, s_{k, \sf r,y}, s_{k, \sf r,z}]$.

In this paper, the users  may need the location information of anchors  to estimate their locations. To preserve the confidentiality of the anchors' locations, they may not share their exact location. According to \cite{dwork2014algorithmic},  a Gaussian noise $n \sim \mathcal{N}(0, \xi_{j,1}^2)$ is added to protect the data of anchor $j$. This noise addition ensures $(\epsilon_j, \delta_j)-$DP if the following condition is met
\begin{equation}
	\begin{aligned}
		\delta_j \geq \frac{4}{5}\exp \left(-\frac{(\xi_{j,1} \epsilon_{j})^2}{2(\Delta_j^{\sf DP})^2}\right),
	\end{aligned}
\end{equation}
where $\Delta_j^{\sf DP}$ represents global sensitivity  of anchor $j$. 
Then, they broadcast their noisy location as
\begin{equation}
	\bs{\check{s}}_{j, \sf a} = \hat{\bs{s}}_{j, \sf a} +\mathcal{N}(0, \xi_{j,1}^2).
\end{equation}
Consequently, the error variance of anchor $j$, denoted as $\xi_j^2 = \xi_{j,0}^2 + \xi_{j,1}^2$.

The anchors can support the location estimation of users.
Let $x_{i,k}$ be a binary indicator where $x_{j,k} = 1$ if the anchor $j$ supports the location estimation of  user $k$. 
Let the set of devices supporting the location estimation of  user $k$ is $\mathcal{X}_k$. Then, $\mathcal{X}_k = \{j\in\mathcal{J}|x_{j,k}=1\}$.
According to \cite{yin2022design}, the positioning error of user $k$ is given by
\begin{equation}
	\Phi_{k} =\sqrt{\sum_{ j \in \mathcal{X}_k} (\lambda_{j,k} \sigma_{j,k})^2},
\end{equation}
where $|.|$ denotes the cardinality of a set, $\lambda_{j,k}$ represents the geometric-dilution and $\sigma_{j,k}$ represents the ranging error variance.
Specifically, let 
\begin{equation}
	\begin{aligned}
		\delta_{j,k}^{\sf x} &= (s_{k, \sf r,x} - s_{j, \sf a,x})/||\bs{s}_{k, \sf r} - \bs{s}_{j, \sf a}||,\\
		\delta_{j,k}^{\sf y} &= (s_{k, \sf r,y} - s_{j, \sf a,y})/||\bs{s}_{k, \sf r} - \bs{s}_{j, \sf a}||,\\
		\delta_{j,k}^{\sf z} &= (s_{k, \sf r,z} - s_{j, \sf a,z})/||\bs{s}_{k, \sf r} - \bs{s}_{j, \sf a}||, \\
			\bs{\check{G}}_{k} &= \begin{bmatrix}
			\delta_{1_{\mathcal{X}_k},k}^{\sf x} & \delta_{1_{\mathcal{X}_k},k}^{\sf y} & \delta_{1_{\mathcal{X}_k},k}^{\sf z} \\ \vdots & \vdots&  \vdots \\ \delta_{|\mathcal{X}_k|_{\mathcal{X}_k},k}^{\sf x} & \delta_{|\mathcal{X}_k|_{\mathcal{X}_k},k}^{\sf y} & \delta_{|\mathcal{X}_k|_{\mathcal{X}_k},k}^{\sf z}
		\end{bmatrix}.
	\end{aligned}
\end{equation}
where $||.||$ represents the Euclidean distance, index $n_{\mathcal{X}_k}$ indicates the $n$th device in set $\mathcal{X}_k$. Then, $\lambda_{j,k}$ is given as follows
\begin{equation}
	\lambda_{j,k} = \sqrt{\sum_{ i=1}^3 g_{i,j,k}^2},
\end{equation}
where $g_{j,i,k}$ represents the element in row $j$ and column $i$ of $\bs{G}_k$. Noted that
\begin{equation}
	\bs{G}_k = \left[\left(\bs{\check{G}}_{k}\right)^T \bs{\check{G}}_{k}\right]^{-1} \left(\bs{\check{G}}_{k}\right)^T.
\end{equation}

In 5G-NR, positioning support involves various comb patterns in the frequency domain, different configurations of symbol numbers in the time domain, and diverse types of numerology. Specifically, 5G NR defines $L=7$ types of numerology, each characterized by a subcarrier spacing (SCS) of $2^l \times 15$ kHz, where $l$ ranges from 0 to 6.  In particular,  let $\omega_j(n,m)$ be the normalized modulated positioning symbol  with subcarrier index $n$ and OFDM symbol index $m$ of anchor $j$. Then, the frequency of the subcarriers carrying normalized symbol $\omega_j(n,m)$ in numerology $l$, denoted as $f_{n,m,l}$, is given by
\begin{equation}
	f_{n,m,l} = (L_{\sf comb} \times n + \kappa_{m,i} )\Delta_{l, \sf f},
\end{equation}
where $\Delta_{l, \sf f} = 1/T_l$ is the subcarrier spacing of numerology $l$, $T_l$  is the duration of the OFDM symbol. $L_{\sf comb}$ is the comb size, and $k_0$ is defined as follows
\begin{equation}
	\kappa_{m,i} = \text{mod}\left(i + \kappa_{m}, L_{\sf comb} \right) - \bar{\kappa}_{m,i},
\end{equation}
where $\kappa_{m} = \frac{\text{mod}\left(m, L_{\sf comb}\right)}{2} + \frac{3}{4}\left(1-\left(-1\right)^{\text{mod}\left(m, L_{\sf comb}\right)}\right)$, mod refers to a modulo operation, $i \in [0, L_{\sf comb}]$, and $\bar{\kappa}_{m,i} = \lfloor \frac{B}{2\Delta_{l,\sf f}} \rfloor$, where $B$ is the available bandwidth. 

In a time slot, the $m$th OFDM symbol occupies $N_l$ subcarriers in the frequency domain when selecting numerology $l$. For each OFDM symbol, the number of subcarriers carrying positioning symbols is $N_{l,\sf a} = N_l/L_{\sf comb}$.  Let $g_{\sf P}(t)$ denote the windowing function.  The transmitted OFDM signal of the $m$th symbol of anchor $j$ when selecting numerology $l$ is then expressed as follows
\begin{equation}
	\begin{aligned}
		s_{j,m,l}(t)=\sum_{n=0}^{N_l-1} c_{j,n,m} g_{\sf P}(t-T) e^{j2\pi n \Delta_{l, \sf f} (t-T)},
	\end{aligned}
\end{equation}
where $c_{j,n,m}$ represents the complex-valued data symbol of anchor $j$ that is modulated on the $n$th subcarrier of the $m$th symbol. For convenience, we denote $\mathcal{N}_{m,i,l} = \{L_{\sf comb} \times n+ \kappa_{m,i}| 0 \leq n \leq N_{l,\sf a}\}$. Let $v_{j,i}$ be a binary variable indicating the selected offset value where $v_{j,i}=1$ if the set of subcarriers carrying normalized $m$th symbol of anchor $j$ is $\mathcal{N}_{m,i,l}$. The mapping condition of positioning symbols on the time-frequency resource domain is expressed as $c_{j,n,m} = 0$ if $n \notin \cap_i v_{j,i} \mathcal{N}_{m,i}$. Only one offset value is selected in the mapping process, which is expressed as follows
\begin{equation}
	\sum_{i=0}^{L_{\sf comb}} v_{j,i} \leq 1, \forall j \in \mathcal{J}. \label{eq12}
\end{equation}

Let $u_{j,l}$ be a binary variable indicating the selected numerology of anchor $j$ where $u_{j,l}=1$ if anchor $j$ selects numerology $l$, otherwise $u_{j,l}=0$. We assume that anchor $j$ can select only one type of numerology, which is expressed as follows
\begin{equation}
	\sum_{l=0}^L u_{j,l} \leq 1, \forall j \in \mathcal{J}. \label{eq13}
\end{equation}


\begin{figure}[t]
	\centering
	\includegraphics[width= 0.98\linewidth]{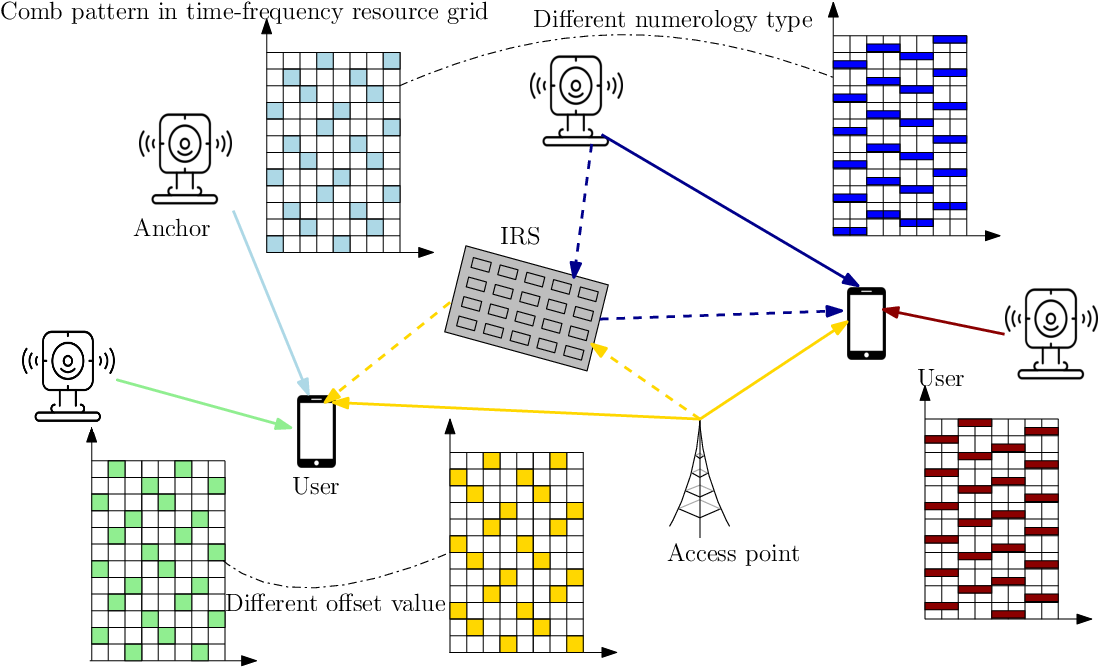}
	\caption{The comb mapping diagram for positioning symbols from various users on the time-frequency resource grid.}
	\label{fig1}
\end{figure}

According to \cite{jayalath2001reducing}, the power spectral density (PSD) of the signal when the comb pattern and rectangular pulse shaping are applied, is given in \eqref{eqlong_1}, where $\sigma_{j,n,m}^2$ is the variance of $c_{j,n,m}$, and $\text{sinc}(x) = \frac{\sin(\pi x)}{\pi x}$.
\begin{figure*}
	\hrulefill
	\begin{equation}
		\begin{aligned}
			P_{j,m}(f) = \sum_{l=0}^{L-1} \sum_{i=0}^{L_{\sf comb}} {T_{l}} u_{j,l}v_{j,i}   \sum_{n \in \mathcal{N}_{m,i,l}} \sigma_{j,n,m}^2 \left(\text{sinc} [(f-n\Delta_{l, \sf f})T_{l}]\right)^2, \label{eqlong_1}
		\end{aligned}
	\end{equation}
\end{figure*}
 Assuming that  $ \sigma_{k,n,m}^2 = \sigma_{k}^2, \forall m, n$ if $c_{j,n,m} \neq 0$. Then the PSD of anchor $j$ on the $m$th symbol is rewritten as
\begin{equation}
	\begin{aligned}
		P_{j,m}(f) \!=\!\! \sum_{l=0}^{L-1}\! T_{l} \sigma_{j}^2\sum_{i=0}^{L_{\sf comb}}  \!  \sum_{n \in \mathcal{N}_{m,i,l}}\!\!\!\!\! u_{j,l} v_{j,i}  \left(\text{sinc} [(f-n\Delta_{l, \sf f})T_{l}]\right)^2.
	\end{aligned}
\end{equation}


%

The received PSD from anchor $j$ to the user $k$ on the $m$th symbol, denoted as $P_{j,k,m}(f)$, is given in \eqref{eqlong_2}, 
\begin{figure*}
	\begin{equation}
		P_{j,k,m}(f) =  \sum_{l=0}^{L-1}  T_{l} \sigma_{j}^2 \sum_{i=0}^{L_{\sf comb}}   \sum_{n \in \mathcal{N}_{m,i,l}} u_{j,l} v_{j,i}   p_{j}|h_{j,k,n}|^2  \left(\text{sinc} [(f-n\Delta_{l, \sf f})T_{l}]\right)^2, \label{eqlong_2}
	\end{equation}
\end{figure*}
where $h_{j,k,n}$ is the channel between anchor $j$ and user $k$ on subcarrier $n$, $p_j$ is the transmit power of anchor $j$.
On the other hands, we have
\begin{equation}
	\begin{aligned}
	\int_{-\infty}^{\infty} P_{j,k,m}(f) 
	&= \sum_{l=0}^{L-1}  {\sigma_{j}^2} \sum_{i=0}^{L_{\sf comb}}   \sum_{n \in \mathcal{N}_{m,i,l}} u_{j,l} v_{j,i}   p_{j}|h_{j,k,n}|^2 \pi. \label{eq16}
		\end{aligned}
\end{equation}
Consequently, its normalized PSD is given by
$\bar{P}_{j,k,m}(f) =   P_{j,k,m}(f)/C_{j,k,m}$, where $C_{j,k,m}$ is a constant such that $\int_{-\infty}^{\infty}\bar{P}_{j,k,m}(f)  df= 1$. From \eqref{eq16}, it can be verify that $C_{j,k,m} = \sum_{l=0}^{L-1}  {\sigma_{j}^2} \sum_{i=0}^{L_{\sf comb}}   \sum_{n \in \mathcal{N}_{m,i,l}} u_{j,l} v_{j,i}   p_{j}|h_{j,k,n}|^2 \pi$. Then, the normalized PSD $\bar{P}_{j,k,m}(f)$ is rewritten as in \eqref{eqlong_3},
\begin{figure*}
	\begin{equation}
		\bar{P}_{j,k,m}(f) =  \bar{\alpha}_{j,k,m}\sum_{l=0}^{L-1}  \frac{T_{l}}{\pi} \sum_{i=0}^{L_{\sf comb}}   \sum_{n \in \mathcal{N}_{m,i,l}} u_{j,l} v_{j,i}   p_{j}|h_{j,k,n}|^2  \left(\text{sinc} [(f-n\Delta_{l, \sf f})T_{l}]\right)^2, \label{eqlong_3}
	\end{equation}
\end{figure*}
where $ \bar{\alpha}_{j,k,m} = \frac{1}{\sum_{l=0}^{L-1}   \sum_{i=0}^{L_{\sf comb}}   \sum_{n \in \mathcal{N}_{m,i,l}} u_{j,l} v_{j,i}   p_{j}|h_{j,k,n}|^2  }$.


The PSD from all transmitters except anchor $j$ to the user $k$   is given in \eqref{eqlong_4},
\begin{figure*}
\begin{equation}
	\begin{aligned}
		\check{P}_{j,k,m}(f) = & \sum_{j'\neq j} \sum_{i =1}^{L_{\sf comb}} \sum_{l=0}^{L-1}  \frac{1}{l_{j,j'}} \sum_{m'=ml_{j,j'}}^{(m+1)l_{j,j'}-1}   v_{j',i} u_{j',l}{T_{l}}{\sigma_{j'}^2} \sum_{n \in \mathcal{N}_{m',i,l}}  p_{j'}|h_{j',k,n}|^2 \left(\text{sinc} [(f-n\Delta_{l, \sf f})T_{l}]\right)^2, \label{eqlong_4}
	\end{aligned}
\end{equation}
\hrulefill
\end{figure*}
where $l_{j,j'} = \max(2^{l-l_j},1)$ such that $u_{j,l_j}=1$. 
According to \cite{yin2022design}, the error variance from anchor $j$ to user $k$  corresponding to slot $m'$ of numerology $0$ is then given in \eqref{eqlong_5},
\begin{figure*}
\begin{equation}
	\begin{aligned}
		\sigma_{j,k}^2 &=  \frac{1}{2^{l_j}} \sum_{m=m'2^{l_j}}^{(m'+1)2^{l_j}-1} \frac{a \int_{-B_{\sf e}/2}^{B_{\sf e}/2} [B\sigma^2 + \check{P}_{j,k,m}(f)] \bar{P}_{j,k,m}(f) \sin^2(\pi f D T) df}{4\pi^2 C_{j,k,m} \left[\int_{-B_{\sf e}/2}^{B_{\sf e}/2} f \bar{P}_{j,k,m}(f) \sin(\pi f D T)df \right]^2}\\
		&= \frac{1}{2^{l_j}} \sum_{m=m'2^{l_j}}^{(m'+1)2^{l_j}-1}  \frac{a(A_{j,k,m,0} + A_{j,k,m,1})}{4\pi^2 C_{j,k,m} A_{j,k,m,2}^2}, \label{eqlong_5}
	\end{aligned}
\end{equation}
\hrulefill
\end{figure*}
where $D$ is the early-late spacing of the delay locked loop (DLL), $B_{\sf e}$ is the double-sided front-end bandwidth, $a = B_{\sf L}(1 - 0.5 B_{\sf L} T_{\sf coh})$ is determined by the loop noise bandwidth $B_{\sf L}$ and the coherent integration time $T_{\sf coh}$, and  
\begin{equation}
	\begin{aligned}
		A_{j,k,m,0} &= \int_{-B_{\sf e}/2}^{B_{\sf e}/2} B\sigma^2 \bar{P}_{j,k,m}(f)  \sin^2(\pi f D T_{l_{j}}) df, \\
		A_{j,k,m,1} &= \int_{-B_{\sf e}/2}^{B_{\sf e}/2}  \check{P}_{j,k,m}(f) \bar{P}_{j,k,m}(f)  \sin^2(\pi f D T_{l_{j}}) df, \\
		A_{j,k,m,2} &=\int_{-B_{\sf e}/2}^{B_{\sf e}/2} f \bar{P}_{j,k,m}(f)  \sin(\pi f D T_{l_{j}})df.
	\end{aligned}
\end{equation}

When $S \rightarrow 0$, $\sin (x) \rightarrow x$, then $A_{j,k,m,i}, i=0,1,2$ are rewritten as in \eqref{eqlong_6},
\begin{figure*}
\begin{equation}
	\begin{aligned}
		A_{j,k,m,0} 
		=& B\sigma^2 \bar{\alpha}_{j,k,m} \sum_{l=0}^{L-1} \pi D^2 T_{l}^3 \sum_{i=0}^{L_{\sf comb}}    \sum_{n \in \mathcal{N}_{m,i,l}} u_{j,l} v_{j,i}  p_{j}|h_{j,k,n}|^2  Q_{n, l},\\
	A_{j,k,m,2} =&  \bar{\alpha}_{j,k,m}\sum_{l=0}^{L-1}  DT_{l}^2 \sum_{i=0}^{L_{\sf comb}}   \sum_{n \in \mathcal{N}_{m,i,l}} u_{j,l} v_{j,i}   p_{j}|h_{j,k,n}|^2 Q_{n, l},\\
		A_{j,k,m,1} =&   \sum_{j'\neq j} \sum_{i =1}^{L_{\sf comb}} \sum_{l=0}^{L-1} \frac{1}{l_{j,j'}} \sum_{m'=ml_{j,j'}}^{(m+1)l_{j,j'}-1}  v_{j',i} u_{j',l} T_{l} {\sigma_{j'}^2} \sum_{n \in \mathcal{N}_{m',i,l}}  p_{j'}|h_{j',k,n}|^2 \\
		& \bar{\alpha}_{j,k,m}\sum_{l_j=0}^{L-1}  {\pi D^2 T_{l_j}^3} \sum_{i_j=0}^{L_{\sf comb}}   \sum_{n_j \in \mathcal{N}_{m,i_j,l_j}} u_{j,l_j} v_{j,i_j}   p_{j}|h_{j,k,n_j}|^2  C_{j',i,l,n}^{j, i_{j},l_{j}, n_{j}}, \label{eqlong_6}
	\end{aligned}
\end{equation}
\hrulefill
\end{figure*}
where $C_{j',i,l,n}^{j, i_{j},l_{j}, n_{j}} = \int_{-B_{\sf e}/2}^{B_{\sf e}/2}  f^2 \left(\text{sinc} [(f-n\Delta_{l, \sf f})T_{l}]\right)^2\left(\text{sinc} [(f-n_{j}\Delta_{l_{j}, \sf f})T_{l_{j}}]\right)^2  df
$, and $Q_{n, l} =  \int_{-B_{\sf e}/2}^{B_{\sf e}/2} f^2 \left(\text{sinc} [(f-n\Delta_{l, \sf f})T_{l}]\right)^2 df$. A detailed analysis of $C_{j',i,l,n}^{j, i_{j},l_{j}, n_{j}}$ and $Q_{n, l}$ is given in Appendix A.


 Let $M_0$ be the number of reflecting elements of IRS, and $\bs{\theta} \in \mathbb{C}^{M_0 \times 1}$ be the coefficients of the reflecting elements of IRS. Assumed $\bs{\theta}$ is generated according to the Kronecker codebook. Then, the channel from the anchor $j$ to user $k$ on the subcarrier $n$ is given by
\begin{equation}
	h_{k,j,n}(t) = h_{k,j,n}^{\sf direct}(t) +  \bs{g}_{k}^n(t) \bs{\Theta}(t) \left(\bs{g}_{j}^n(t) \right)^H,
\end{equation}
where $\bs{\Theta}(t) = \text{diag}(\bs{\theta}) \in \mathbb{C}^{M_0 \times M_0}$, $\bs{g}_{k}^n(t)$ are the channel between IRS and user/anchor $k$ on subcarrier $n$, and $h_{k,j,n}^{\sf direct}(t)$ is the direct channel between anchor $j$ to user $k$ on the subcarrier $n$.

Aiming to achieve fairness for different users, we formulate the problem that minimizes the maximum users' positioning error in an IRS-assisted 5G network as follows
\begin{subequations}
	\begin{alignat}{4}
		(\mathcal{P}) &    \min\limits_{\bs{X}, \bs{U}, \bs{V},\bs{\theta},\bs{p}}  \; \;  \max_{k \in \mathcal{K}} \Phi_{k}(t) 	\nonumber \\
		\text{ s.t. }  
		& \bs{\theta} \in \mathcal{C}, \label{eqbeam}\\
		& \sum_{l=0}^L u_{j,l} N_{l,a} p_j \leq P_j^{\max}, \forall j \in \mathcal{J}, \label{eqpower}\\
		& \xi_j^2 \geq \xi_{j,\min}^2, \forall j \in \mathcal{J}, \label{eqlong_12}\\
		& \eqref{eq12}, \eqref{eq13}, \nonumber
	\end{alignat}
\end{subequations} 
where $\bs{X} = \{x_{j,k}, \forall j,k\}$, $\bs{U} = \{u_{j,l}, \forall j,l\}$, $\bs{V} = \{v_{j,i}, \forall j,i\}$, and $\bs{p} = \{p_{k}, \forall k\}$, $\mathcal{C}$ is the set of Kronecker codewords, $P_j^{\max}$ is the maximum transmit power of anchor $j$, and $\xi_{j,\min}^2$ is the minimum anchor error variance which is determined by $\xi_{j,0}^2 - 2(\Delta_j^{\sf DP})^2/\epsilon_{j,\min}^2 \log \left( 5\delta_{j,\min}/4 \right)$, $(\epsilon_{j,\min}, \delta_{j,\min})$ are the minimum value of $(\epsilon_{j}, \delta_{j})$-DP. In problem $(\mathcal{P})$, \eqref{eqbeam} indicates the IRS coefficient vector is selected from the predefined codebook $\mathcal{C}$. The anchors' maximum power constraint is presented in \eqref{eqpower}, while \eqref{eqlong_12} captures the anchors'  privacy condition.


Due to binary variables $u_{j,l}$, $v_{j,i}, \forall j,i,l$, constraints \eqref{eqbeam}, \eqref{eqpower}, the non-convex  function $\Phi_{k}(t)$,  and the min-max objective, $(\mathcal{P}) $ is a MINLP problem. In the next section, we present an efficient algorithm to tackle this challenging problem.

\section{Hybrid optimization,  matching, and DQN approach (HOMD)}
\label{sec_4}

In this section, we introduce a low-complexity algorithm, named hybrid optimization,  matching, and DQN (HOMD), which iteratively solves four subproblems decomposed from problem $\mathcal{P}$ until convergence. Specifically, the HOMD mechanism employs a DQN-based method for the beam selection according to which the transmission power and privacy protection level are optimized efficiently. Additionally, numerology and offset value selection are performed through a many-to-one matching process, while user-anchor associations are established using a many-to-many matching technique.

\subsection{Numerology and offset value selection}
\label{sec_4a}
Let $\mathcal{L}$ and $ \mathcal{I}$ be the set of numerology types and offset values, respectively.
Let define the combined set $\mathcal{D} = \{ l \in \mathcal{L}, i \in \mathcal{I}\}$.

\begin{definition}
	A many-to-one matching $\mu$ is a subset in $ \mathcal{J} \times \mathcal{D}$ such that $|\mu(j)|\leq 1$ and $|\mu((l,i))| \leq J$, where $\mu(j) = \{(l,i) \in \mathcal{D}: (j,(l,i)) \in \mu\}$ and $\mu((l,i)) = \{j \in \mathcal{J}: (j,(l,i)) \in \mu\}$.
\end{definition}

According to \cite{zhao2016many}, the proposed matching game has externalities, where the preference value for anchor $j$ not only depends on its selection that it is matched with but also on the matching of other $j' \in \mathcal{J}$.
We define the preference value for the anchor $j$ selecting $(l,i) \in \mathcal{D}$ as 
$
		U_{j}^{\sf o}( \mu) =   \max_{k \in \mathcal{K}_{j}}  \sqrt{\sum_{ j \in \mathcal{X}_k} (\lambda_{j,k} \sigma_{j,k})^2},
$
where $\mathcal{K}_j = \{k \in \mathcal{K}|j \in  \mathcal{X}_k \}$.
To find the solution of the matching game with externalities,
 we define the preference value $(l,i) \in \mathcal{D}$  selected by anchors as 
 \begin{equation}
 	\begin{aligned}
 		U_{ (l,i)}^{\sf o}(\mu) =   \max_{k \in \mathcal{K}}  \sqrt{\sum_{ j \in \mathcal{X}_k} (\lambda_{j,k} \sigma_{j,k})^2}.
 	\end{aligned}
 \end{equation}

To define the exchange stability, we first consider the concept of swap matching as follows:
\begin{equation}
	\mu_{j}^{j'} = \mu  \backslash \{(j,\mu(j)), (j', \mu(j'))\} \cup \{(j,\mu(j')), (j', \mu(j))\}.
\end{equation}
In other words, a swap matching enables switching from $(j,\mu(j))$ to $(j,\mu(j'))$ and from $(j',\mu(j'))$ to $(j',\mu(j))$ while keeps all other matchings. Further, the anchor involved in the swap can be a hole `$\mathcal{O}$' representing an open spot, thus allowing a single $(l,i) \in \mathcal{D}$ moving to available vacancies.

\begin{definition}
	Let $j \in \mathcal{J}$ and $j' \in \mathcal{J}\cup \mathcal{O}$, $(j, j')$ is a swap-blocking pair if and only if 1) $\forall x \in \{ j, (l,i) = \mu(j), j', (l',i') = \mu(j')\} $, $U_{x}^{\sf o}(\mu_{j}^{j'}) \leq U_{x}^{\sf o}(\mu)$, and 2) $\exists x \in \{ j, (l,i) = \mu(j), j', (l',i') = \mu(j')\} $ such that $U_{x}^{\sf o}(\mu_{j}^{j'}) < U_{x}^{\sf o}(\mu)$.
\end{definition}

According to \cite{zhao2016many}, a matching $\mu$ is a two-sided exchange-stable if no swap-blocking pairs exists. Therefore, we propose the algorithm to determine the binary decision $u_{j,l}$ and $v_{j,i}, \forall j \in \mathcal{J} $ that attains the two-sided exchange-stable matching as in Algorithm~\ref{alg2}. Indeed, if the converged matching $\mu$ is not the two-sided exchange-stable matching, a swap-blocking pair exists after the final matching $\mu$. However, as shown in Algorithm~\ref{alg2}, the algorithm does not terminate until all swap-blocking pairs are eliminated. Therefore, the converged $\mu$ must be a two-sided exchange-stable matching.

Furthermore, the decrease of positioning error after each swap operation shows the convergence of Algorithm~\ref{alg2}. According to Algorithm~\ref{alg2}, a swap occurs only when $U_{(l,i)}(\mu_{j}^{j'}) \leq U_{(l,i)}(\mu)$. Then, \begin{equation}
	\Phi_{\mu \rightarrow \mu_{j}^{j'}} {=} \max_k \Phi_{k}|_{\mu_{j}^{j'}} - \max_k \Phi_{k}|_{\mu} {=} U_{ (l,i)}^{\sf o}(\mu_{j}^{j'}) - U_{ (l,i)}^{\sf o}(\mu) \leq 0,
\end{equation}
 In this case, $\Phi_{\mu \rightarrow \mu_{j}^{j'}}$ is the difference of positioning error under the matching
state $\mu_{j}^{j'}$ and that under the matching state $\mu$

\setlength{\textfloatsep}{5 pt}
\begin{algorithm}[t]
	\footnotesize
	\caption{Numerology and offset selection algorithm }
	\label{alg2}
	\begin{algorithmic}[1]
		\State \textbf{Initialize}: Variable $u_{j,l}, v_{j,i}, \forall j,l,i$ are randomly assigned.
		\State \textbf{repeat}
		\State For each $j$, search for another $j'$ or an open spot $\mathcal{O}$ to form a swap-blocking pair.

		\If{$\{j,j'\}$ or $\{j,\mathcal{O}\}$ forms a swap-blocking pair}
		\State Update $\mu = \mu_{j}^{j'}$ or $\mu = \mu_{j}^{\mathcal{O}}$, respectively.
		\EndIf 
		
		\State \textbf{until} matching $\mu$ is two-sided exchange-stable.
	\end{algorithmic}
\end{algorithm}

\setlength{\textfloatsep}{5 pt}
\begin{algorithm}[t]
	\footnotesize
	\caption{User-anchor association algorithm }
	\label{alg3}
	\begin{algorithmic}[1]
		\State \textbf{Initialize}: Variable $x_{j,k}$ are randomly assigned such that $\sum_{j \in \mathcal{J}} x_{j,k} \geq 3, \forall k \in \mathcal{K}$.
		\State \textbf{repeat}
		\State For each $j$, search for another $j'$ or an open spot $\mathcal{O}$ to form a swap-blocking pair.

		\If{$\{j,j'\}$ or $\{j,\mathcal{O}\}$ forms a swap-blocking pair}
		\State Update $\nu = \nu_{jk}^{j'k'}$ or $\nu = \nu_{jk}^{\mathcal{O}k'}$, respectively.
		\EndIf 
		
		\State \textbf{until} matching $\nu$ is two-sided exchange-stable.
	\end{algorithmic}
\end{algorithm}
\subsection{Determination of user-anchor association}
\label{sec_4b}
\begin{definition}
	A many-to-many matching $\nu$ is a subset in $\mathcal{J} \times \mathcal{K}$ such that $ |\nu(j)|\leq K$ and $3 \leq |\nu(k)| \leq J$, where $\nu(j) = \{k \in \mathcal{K}: (j,k) \in \nu\}$ and $\nu(k) = \{j \in \mathcal{J}: (j,k) \in \nu\}$.
\end{definition}

This proposed matching game has externalities.
We define the preference value for user $k$ selecting anchors $j \in \mathcal{K}$ as 
$
		U_{k}^{\sf A}( \nu) =  \sum_{ j \in \mathcal{X}_k}  \lambda_{j,k} \sigma_{j,k}.
$
To find the solution of the matching game with externalities,
we define the preference value for anchor $j \in \mathcal{J}$  selected by users $k \in \mathcal{K}_{\nu}$ as 
\begin{equation}
	\begin{aligned}
		U_{j}^{\sf A}(\nu) =   \max_{k \in \mathcal{K}_{\nu}}  \sqrt{\sum_{ j \in \mathcal{X}_k} (\lambda_{j,k} \sigma_{j,k})^2}.
	\end{aligned}
\end{equation}

To define the exchange stability, we first consider the concept of swap matching as follows:
\begin{equation}
	\nu_{jk}^{j'k'} = \nu  \backslash \{(j,k), (j', k')\} \cup \{(j,k'), (j', k)\},
\end{equation}
where $(j,k) \in \nu$ and $(j',k') \in \nu$.
In other words, a swap matching enables switching from $(j,k)$ to $(j,k')$ and from $(j',k')$ to $(j',k)$ while keeps all other matchings. Further, the anchor involved in the swap can be a hole `$\mathcal{O}$' representing an open spot, thus allowing a single user $k \in \mathcal{K}$ moving to available vacancies.

\begin{definition}
	Let $j \in \mathcal{J}$ and $j' \in \mathcal{J}\cup \mathcal{O}$, $(j, j')$ is a swap-blocking pair if and only if 1) $\forall x \in \{ j, k, j', k'\} $, $(j,k), (j',k') \in \nu$, $U_{x}^{\sf A}(\nu_{jk}^{j'k'}) \leq U_{x}^{\sf A}(\nu)$, and 2) $\exists x \in \{ j, k, j', k'\} $, $(j,k), (j',k') \in \nu$ such that $U_{x}^{\sf A}(\nu_{jk}^{j'k'}) < U_{x}^{\sf A}(\nu)$.
\end{definition}

Similarity, a matching $\nu$ is also a two-sided exchange-stable if no swap-blocking pairs exists. Therefore, we propose the algorithm to determine the binary decision $x_{j,k}$ that attains the two-sided exchange-stable matching as in Algorithm~\ref{alg3}.

\subsection{Determination of transmit power and privacy protection level}
\label{sec_4c}
For given $\bs{\theta}$, $\bs{x}$, $u_{j,l_k^{\star}}=1$ and $v_{j,i_k^{\star}}=1, \forall j \in \mathcal{J}$,  the  error variance from anchor $j$ to user $k$  is rewritten as in \eqref{eqlong_20},
\begin{figure*}
	\hrulefill
	\begin{equation}
		\begin{aligned}
			\sigma_{j,k}^2 &=  \frac{1}{2^{l_j}} \sum_{m=m'2^{l_j}}^{(m'+1)2^{l_j}-1} \frac{\zeta_{j,k, l_j^{\star},i_j^{\star},m,0} + \sum_{j'\neq j}  \frac{1}{l_{j,j'}} \sum_{m'=ml_{j,j'}}^{(m+1)l_{j,j'}-1}  p_{j'}   \zeta_{j,j', k, l_j^{\star},i_j^{\star},l_{j'}^{\star},i_{j'}^{\star},m',m}}{  p_j \zeta_{j,k, l_j^{\star},i_j^{\star},m,2}^2}, \label{eqlong_20}
		\end{aligned}
	\end{equation}
	\hrulefill
\end{figure*}
where
\begin{equation}
	\begin{aligned}
		\zeta_{j,k, l_j^{\star},i_j^{\star},m,0}
		&=  a B\sigma^2  \pi   \sum_{n \in \mathcal{N}_{m,i_j^{\star},l_j^{\star}}}   |h_{j,k,n}|^2  Q_{n, l_j^{\star},i_j^{\star}}, \\
	\zeta_{j,k, l_j^{\star},i_j^{\star},m,2} &= \frac{\sqrt{4\sigma_j^2 \pi^3}}{T_{l_j}}   \sum_{n \in \mathcal{N}_{m,i_j,l_j}}   |h_{j,k,n}|^2 Q_{n, l_j, i_j},	
	\end{aligned}
\end{equation}
\begin{equation}
	\begin{aligned}
		&\zeta_{j,j', k, l_j^{\star},i_j^{\star},l_{j'}^{\star},i_{j'}^{\star},m',m} = a {\pi T_{l_j}^3} T_{l_{j'}} {\sigma_{j'}^2} \\
		&   \hspace{0.25 cm}  \sum_{n \in \mathcal{N}_{m',i_{j'},l_{j'}}}         \sum_{n_j \in \mathcal{N}_{m,i_j,l_j}} |h_{j',k,n}|^2 |h_{j,k,n_j}|^2  C_{j',i,l,n}^{j, i_{j},l_{j}, n_{j}}.
	\end{aligned}
\end{equation}
Then, for given $\bs{\theta}$, $\bs{x}$, $u_{j,l_k^{\star}}=1$ and $v_{j,i_k^{\star}}=1, \forall j \in \mathcal{J}$, the problem is written as follows
\begin{subequations}
	\begin{alignat}{4}
		(\mathcal{P}) &    \min\limits_{ \bs{\xi}^2, \bs{p}}  \; \;  \max_k \sum_{ j \in \mathcal{X}_k} \lambda_{j,k}^2 (\xi_j^2 + \sigma_{j,k}^2)
		\text{ s.t. }  
		\eqref{eqlong_12}. \nonumber
	\end{alignat}
\end{subequations} 

In order to convexify the non-convex objective, we put $p_j = \exp(\tilde{p}_j)$. Then, the problem is equivalent rewritten as follows
\begin{subequations}
	\begin{alignat}{4}
		(\mathcal{P}_1) &    \min\limits_{ \bs{\xi}^2, \tilde{\bs{p}}, \varepsilon}  \varepsilon	\nonumber \\
		\text{ s.t. }  
		&\eqref{eqlong_12}; \eqref{eqlong_11},  \forall k \in \mathcal{K}. \nonumber 
	\end{alignat}
\end{subequations} 

\begin{figure*}
	\begin{equation}
		\sum_{ j \in \mathcal{X}_k} \lambda_{j,k}^2 \left(\xi_j^2 +  \frac{1}{2^{l_j}} \sum_{m=m'2^{l_j}}^{(m'+1)2^{l_j}-1} \frac{\exp(-\tilde{p}_j)\zeta_{j,k, l_j^{\star},i_j^{\star},m,0} + \sum\limits_{j'\neq j}  \frac{1}{l_{j,j'}} \sum_{m'=ml_{j,j'}}^{(m+1)l_{j,j'}-1}  \exp(\tilde{p}_{j'}-\tilde{p}_j)   \zeta_{j,j', k, l_j^{\star},i_j^{\star},l_{j'}^{\star},i_{j'}^{\star},m}}{   \zeta_{j,k, l_j^{\star},i_j^{\star},m,2}^2} \right) \leq \varepsilon, \label{eqlong_11}
	\end{equation}
\hrulefill
\end{figure*}
According to \cite{boyd2004convex}, problem $(\mathcal{P}_1)$ is a standard convex problem which can effectively solved by the CVX-solver.

\subsection{IRS's beam selection}
\label{sec_4d}
Employment DQN facilitates efficient learning of beam decisions. Particularly, the state, action, and reward is defined as follows. The state $\bs{s}(t)$ includes the information about channel gain between users with the beam obtained from the action in previous time-slot. Specifically, let
\begin{equation}
	\begin{aligned}
	h_{k,j,n}^{\sf delay}(t) &= h_{k,j,n}^{\sf direct}(t) +  \bs{g}_{k}^n(t) \bs{\theta}(t-1) \left(\bs{g}_{j}^n(t) \right)^H,\\
	\bs{h}_{k,j}^{\sf delay}(t)&=[h_{k,j,1}^{\sf delay}(t),h_{k,j,2}^{\sf delay}(t),..., h_{k,j,N}^{\sf delay}(t)]
	\end{aligned}
\end{equation}
Then, state $\bs{s}(t)$ is given by
\begin{equation}
	\bs{s}(t) = [||\bs{h}_{k,j}^{\sf delay}(t)||_2], \forall k \in \mathcal{K}, j \in \mathcal{J}.
\end{equation}

Additionally, one defines the action of the learner for the beam selection, expressed as
\begin{equation}
	a(t) = \{1,2,..., N_{\sf cb}\},
\end{equation}
where $N_{\sf cb}$ is the size of the codebook $\mathcal{C}$.
The reward is the positioning accuracy $r(t) = 1-\max_k \Phi_k$.

Based on the actions and rewards obtained from trials, the DQN agent builds its DQN model consisting of two deep neural networks (DNNs), namely online and target networks corresponding to weight vectors $\bs{\phi}$ and $\bs{\phi}'$, respectively. Herein, the online network is used to select an action. Meanwhile, the target network is applied to evaluate the online network-based action. Thus, the  objective is to reduce the loss function as 
\begin{equation}
	\begin{aligned}
		L(\bs{\phi}) = \left( y(t) - Q(\bs{s}(t), a(t); \bs{\phi})\right)^2
	\end{aligned}
\end{equation}
where $y(t)$ denotes the target Q-value determined by the target network as $y(t) = r(t) + \max_{a \in \mathcal{A}} Q(\bs{s}(t+1), a; \bs{\phi}')$. The network has 3 hidden layers with size of 1024, 512, and 512, respectively. The ReLu activation is applied at hidden layers, while the sigmoid activation is applied at the output layer. The Adam algorithm is employed in the optimizer with learning rate is set to 0.001.

Finally, the HOMD approach is summarized in Algorithm~\ref{alg4}.

\setlength{\textfloatsep}{5 pt}
\begin{algorithm}[t]
	\footnotesize
	\caption{Hybrid optimization, matching, and DQN approach (HOMD) }
	\label{alg4}
	\begin{algorithmic}[1]
		\State \textbf{Initialize}: Variable $\bs{X}, \bs{U}, \bs{V}$ are randomly assigned, $\bs{p} = \bs{p}^{\max}, \bs{\sigma} = - \frac{2\Delta^2}{\epsilon^2} \log \frac{5\delta}{4}$.
		\State \textbf{repeat}
		\State \hspace{0.25 cm} Solve $\mathcal{P}_1$ to get the transmit power and privacy protection level as presented in Section~\ref{sec_4c}.
		\State  \hspace{0.25 cm}  Determine the user-anchor association as presented in Section~\ref{sec_4b}.
		\State  \hspace{0.25 cm}  Determine the numerology and offset values of anchors as presented in Section~\ref{sec_4a}.
		\State  \hspace{0.25 cm}  Select the beam of IRS as presented in Section~\ref{sec_4d}.
		\State \textbf{until} convergence.
	\end{algorithmic}
\end{algorithm}
\section{Numerical results}

\label{sec_5}
In this section, we numerically evaluate the performance of the proposed framework. Consider a scenario with $J$ anchors and $K$ users uniformly distributed at random within a cell coverage area measuring $100 \times 100$ square meters. The access point (AP) is located at coordinates $(30, 30, 3)$ meters, while the IRS is positioned at (50,50,3) meters within the same area. 
The pathloss (PL) between the user and anchor is defined as the combination of line-of-sight (LoS) and non-line-of-sight (NLoS) components, represented as $\text{PL} = \text{Pr}_{\sf LoS} \text{PL}_{\sf LoS} + (1-\text{Pr}_{\sf LoS}) \text{PL}_{\sf NLoS} $, where $\text{Pr}_{\sf LoS}$ is the line-of-sight (LoS) probability, $\text{PL}_{\sf LoS}$ is the LoS pathloss, and $\text{PL}_{\sf NLoS}$ is the non-LoS pathloss. According to ETSI TR 138901 \cite{3GPP_channel}, we consider the `mixed mode' and `open mode' scenarios to model the LoS probability. In particular, let $d_{\sf 2D\_in}$  and $d_{\sf 3D}$ be the 2D\_in and 3D distances\cite{3GPP_channel} between the transmitter and receiver, respectively; the LoS probability in `mixed mode', denoted as $\text{Pr}_{\sf LoS}^{\sf MM}$, is given by
\begin{equation}
	\begin{aligned}
		\text{Pr}_{\sf LoS}^{\sf MM} = \begin{cases}
			1, \text{ if } d_{\sf 2D\_in} \leq 1.2 \text{m} \\
			\exp\left(-\frac{d_{\sf 2D\_in}-1.2}{4.7}\right), \text{ if } 1.2 \text{m} \leq d_{\sf 2D\_in} \leq 6.5 \text{m} \\
			\exp\left(-\frac{d_{\sf 2D\_in}-6.5}{32.6}\right)0.32, \text{ if } 6.5 \text{m} \leq d_{\sf 2D\_in}. 
		\end{cases}
	\end{aligned}
\end{equation}
The LoS probability in `open mode', denoted as $\text{Pr}_{\sf LoS}^{\sf OM}$, is given by
\begin{equation}
	\begin{aligned}
		\text{Pr}_{\sf LoS}^{\sf OM} = \begin{cases}
			1, \text{ if } d_{\sf 2D\_in} \leq 5 \text{m} \\
			\exp\left(-\frac{d_{\sf 2D\_in}-5}{70.8}\right), \text{ if } 5 \text{m} \leq d_{\sf 2D\_in} \leq 49 \text{m} \\
			\exp\left(-\frac{d_{\sf 2D\_in}-49}{211.7}\right)0.54, \text{ if } 49 \text{m} \leq d_{\sf 2D\_in}. 
		\end{cases}
	\end{aligned}
\end{equation}
The LoS and NLoS pathloss are respectively calculated as follows
\begin{equation}
	\begin{aligned}
		\text{PL}_{\sf LoS}\! &{=} \!32.4 + 17.3\log_{10}(d_{\sf 3D})  +20\log_{10}(f_{\sf c}), \\
		\text{PL}_{\sf NLoS}\! &{=} \!\max \left(\text{PL}_{\sf LoS}, 17.3 {+} 38.3\log_{10}(d_{\sf 3D})  {+}24.9\log_{10}(f_{\sf c})\right), \nonumber
	\end{aligned}
\end{equation}
where $f_{\sf c}$ is the center frequency in Hz.

The channel between user $k$ and anchor $j$ is modeled as $h_{j,k,n}^{\sf direct} = u_{j,k,n}\sqrt{\beta_{j,k}} $ where $ u_{j,k,n} \sim \mathcal{CN}(0,1)$, and $\beta_{j,k} = 10^{-\text{PL}_{j,k}/10}$, where $\text{PL}_{j,k}$ is the pathloss between user $k$ and anchor $j$.
The channel between IRS and user/anchor $k$ is
$
\bs{g}_{k}^n = \bs{u}_{k,n}^{\sf IRS} (\psi_{k,1}, \psi_{k,2}) \sqrt{\beta_{k, \sf LoS}},
$
where $(\psi_{k,1}, \psi_{k,2}) $ are the (azimuth, elevation) angles, $\bs{u}_{k,n}^{\sf IRS}$ is the steering vector and $\beta_{k, \sf LoS}  = 10^{-\text{PL}_{k,\sf LoS}/10}$, where $\text{PL}_{k,\sf LoS}$ is the LoS pathloss between the IRS and user $k$.

\begin{figure}[t]
	\centering
	\includegraphics[width= 0.9 \linewidth]{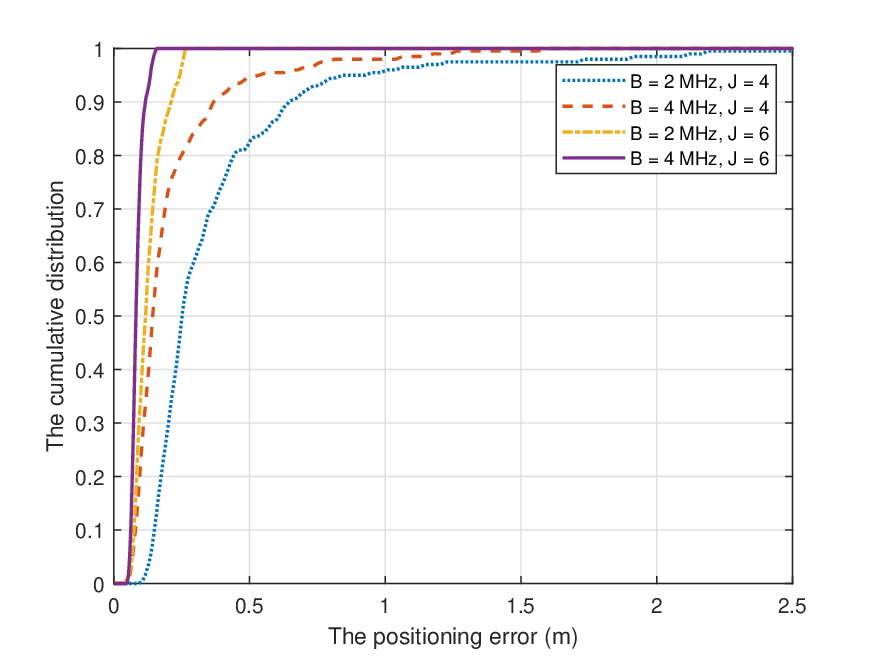}
	\caption{The cumulative distribution of the positioning error.}
	\label{JFig5}
\end{figure}

In all experiments, we set the default values as follows: the number of users $K = 9$, the number of anchors $J = 6$, `open mode' for LoS probability between users and anchors, the bandwidth $B = 4$ MHz, the maximum transmit power of the AP and users, denoted as $\bs{p}_0^{\max}$, are 0.22 Watts and 0.05 Watts, respectively. In default, we set $\bs{P}^{\max} = \bs{p}_0^{\max}$. The IRSs' antennas is $M_0 = 5 \times 5$. 
The noise power is set as
$
\text{noise power} = B \times 3.1812 \times 10^{-20}
$. The minimum anchor error variance is $\xi_{j,\min}^2 = 0.005$ m. The signal is modulated by 16-QAM technology. 
The loop parameters are set as
$D = 0.02$ chips, $B_{\sf e} = 2 B$ , $a = B_{\sf L}(1 - 0.5 B_{\sf L} T_{\sf coh})$,  $B_{\sf L} = 0.2$ Hz, and $T_{\sf coh} = 0.02$ s. The simulation results are obtained by averaging the results over 100 realizations except for Fig.~\ref{JFig5}.

\begin{figure}[t]
	\centering
	\includegraphics[width= 0.9 \linewidth]{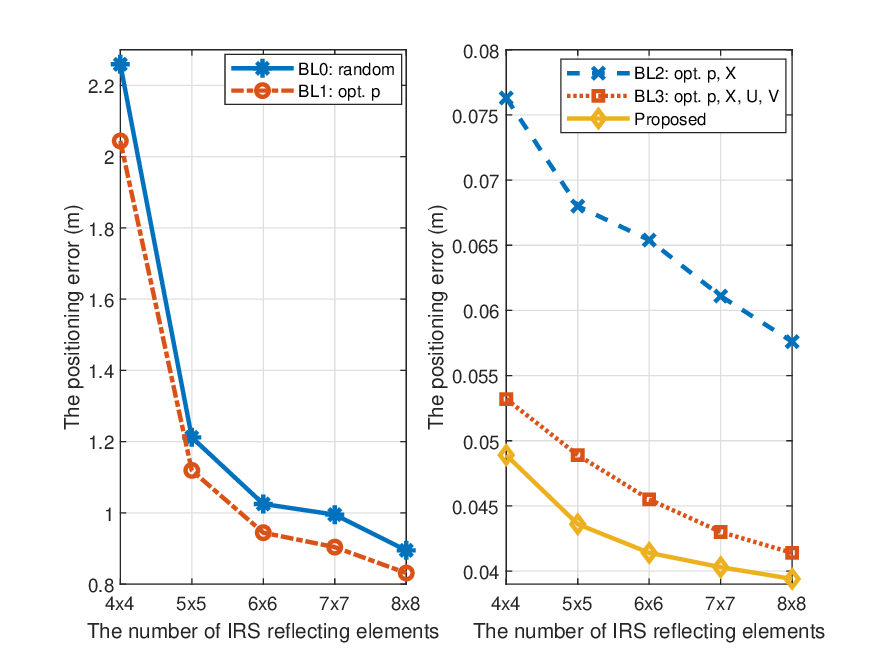}
	\caption{The impact of the number of reflecting elements on the positioning error.}
	\label{JFig2}
\end{figure}

In Fig.~\ref{JFig2},  the impact of the number of IRS reflecting elements on the positioning error is demonstrated. To evaluate the proposed frameworks, we also plot the positioning error obtained in four baseline (BL) schemes: `BL0: random', `BL1: opt. p', `BL2: opt. p, X', and `BL3: opt. p, X, U, V'. Specifically, `BL0: random' represents results obtained when numerology, offset value, and user-anchor associations are randomly chosen, with the condition that each user is associated with at least 3 anchors. Additionally, the curve incorporates the application of maximum transmit power and minimum privacy protection levels. In `BL1: opt. p', all variables remain the same as in the `BL0: random' case, except the transmit power and privacy protection level, which are determined by the algorithm proposed in Section~\ref{sec_4c}. In `BL2: opt. p, X', all parameters are retained from the `BL1: opt. p' case, except the user-anchor association, which is determined using the algorithm outlined in Section~\ref{sec_4b}. Further refinement occurs in the `BL3: opt. p, X, U, V' scenario, where all variables remain the same as in `BL2: opt. p, X', except the numerology and offset value selection, which are optimized using the algorithm detailed in Section \ref{sec_4a}. In `Proposed', all variables are determined by our proposed HOMD algorithm presented in Section~\ref{sec_4}. Fig.~\ref{JFig2} underscores the significance of controlling these variables to minimize positioning errors. Moreover, the results obtained through the proposed optimization, matching, and learning algorithm are up to 25 times superior to those obtained without the application of these proposed techniques.

\begin{figure}[t]
	\centering
	\includegraphics[width= 0.9 \linewidth]{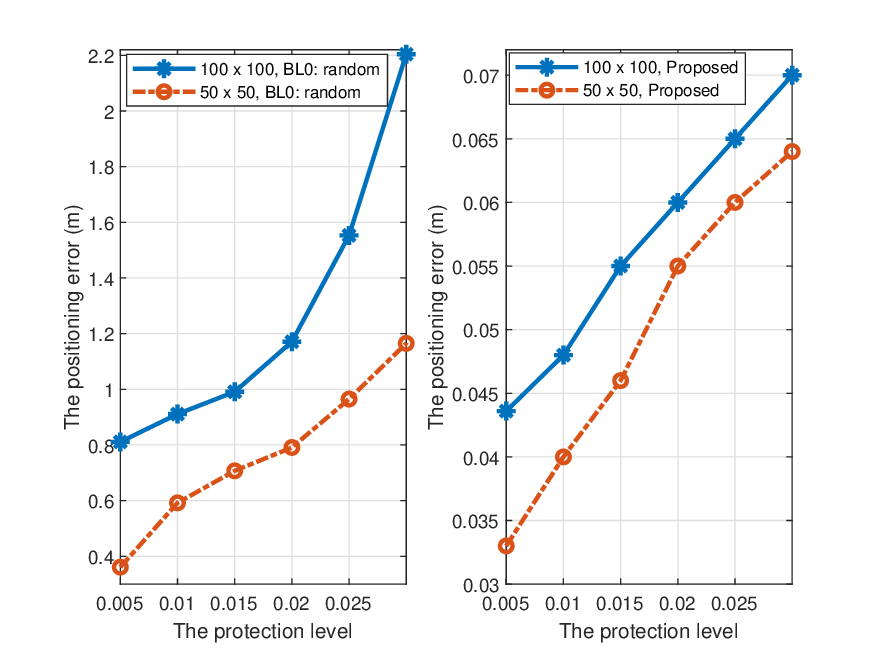}
	\caption{The impact of differential privacy protection level on the positioning error.}
	\label{JFig3}
\end{figure}
In Figure \ref{JFig3}, we illustrate the influence of the differential privacy protection level on positioning errors within both a $100 \times 100$ $\text{m}^2$ and a $50 \times 50$ $\text{m}^2$ area. As the differential privacy protection level rises, there is a corresponding increase in positioning errors.  Additionally, when the distance decreases, there is a higher likelihood of LoS conditions, resulting in increased channel gains and consequently, reduced positioning errors.

\begin{figure}[t]
	\centering
	\includegraphics[width= 0.9 \linewidth]{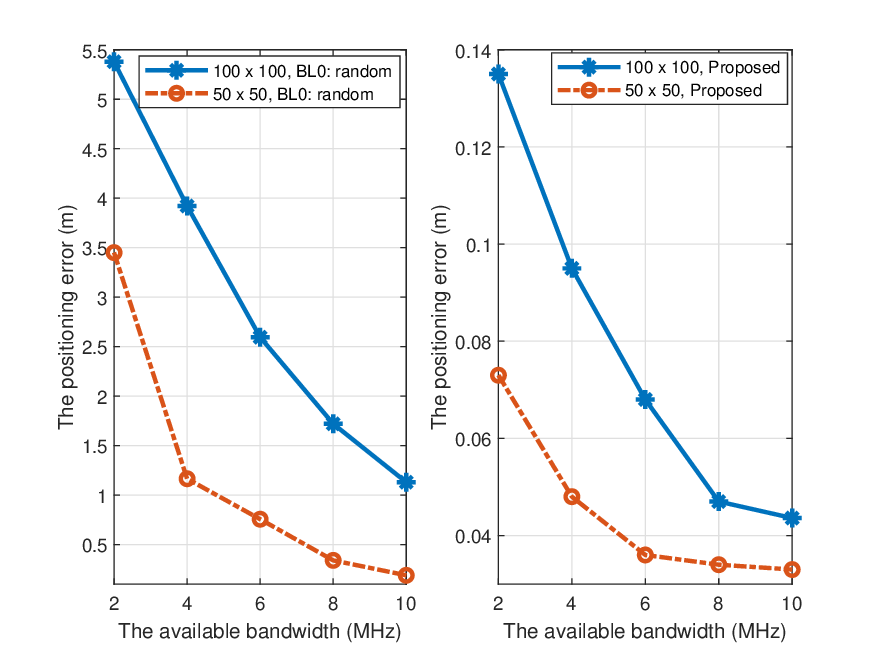}
	\caption{The impact of bandwidth on the positioning error.}
	\label{JFig4}
\end{figure}

In Fig.~\ref{JFig4}, the influence of bandwidth on positioning error is depicted. Notably, when the bandwidth is limited to $2$ MHz, the positioning accuracy is achieved at the decimeter level. However, with an increase in bandwidth to $10$ MHz, the positioning estimation markedly improves, approaching accuracy levels close to the millimeter scale.

Fig~\ref{JFig5}  presents the cumulative distribution of positioning errors under varying numbers of anchors and available bandwidth. In scenarios where both the bandwidth and the number of anchors are limited, location estimations can be imprecise, with errors extending up to 2.5 meters within a $100 \times 100 , \text{m}^2$ area. However, as the bandwidth and the number of anchors increase, the scenario significantly improves. In cases where $B = 4$ MHz and $J=6$, more than 90\% of positioning estimations achieve an error of less than 1 decimeter. The results in Fig~\ref{JFig5} are obtained over 1000 realizations.

\begin{figure}[t]
	\centering
	\includegraphics[width= 0.9 \linewidth]{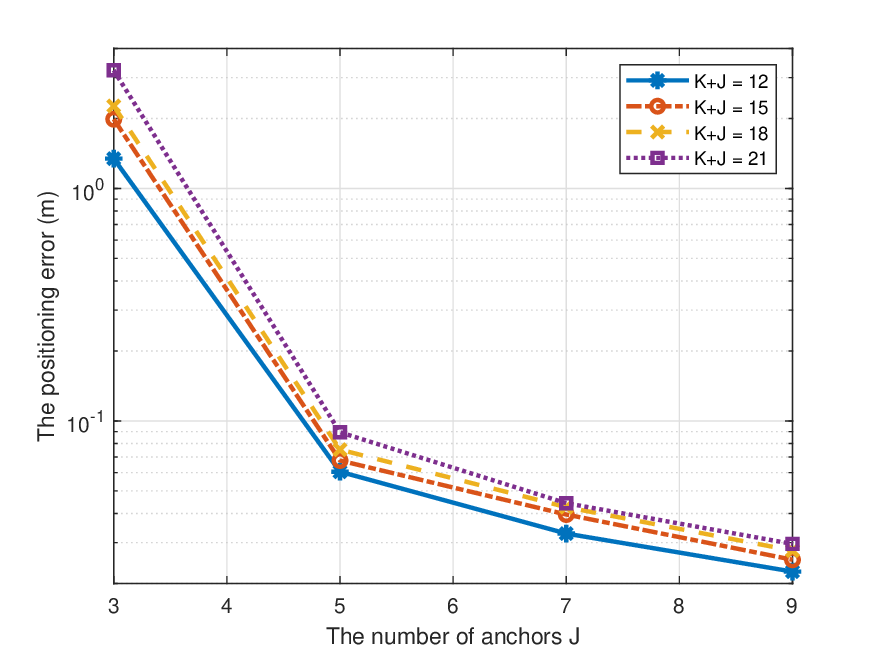}
	\caption{The positioning error with different number of anchors and users.}
	\label{JFig6}
\end{figure}

In Fig~\ref{JFig6}, the positioning errors are depicted under different number of anchors and users configurations. Specifically, when the number of anchors is limited, such as  $J=3$, doubling the number of users (from 9 to 18) results in a twofold increase in positioning error if the wireless resources remain unchanged. When the number of anchor is larger (i.e, $J = 9$), increasing the number of users fourfold (from $3$ to $12$) only leads to a decrease in accuracy of less than 30\%.

\begin{figure}[t]
	\centering
	\includegraphics[width= 0.9 \linewidth]{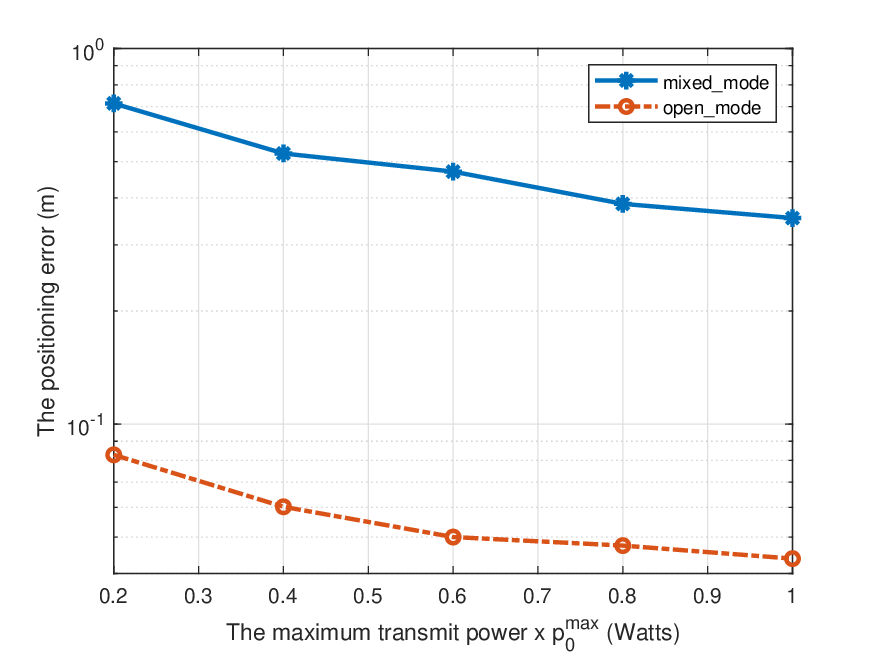}
	\caption{The positioning error in different environments.}
	\label{JFig8}
\end{figure}

In Fig.~\ref{JFig8}, the positioning error is depicted for scenarios where the channel LoS probability model alternates between `mixed mode' and `open mode'. In the `mixed mode', the LoS probability is lower than that in the `open mode'. Consequently, in the `mixed mode', a positioning accuracy at the decimeter level is achieved when the maximum transmit power $p_0^{\max}$ is employed. Conversely, in the `open mode', characterized by a higher LoS probability, the positioning error can be reduced to the centimeter level.
\section{Conclusions}
\label{sec_6}
In this paper, we introduce an innovative design for multi-user positioning in 5G NR. Our approach addresses several critical aspects, including joint user-anchor association, offloading and offset value decision, transmit power and privacy protection level assignment, as well as beam selection to control IRS elements. The objective is to minimize the maximum users' positioning error within a 5G system. Our numerical results demonstrate the effectiveness of our proposed solution, showcasing a substantial reduction in positioning error within a 5G network.

\bibfont
\footnotesize{
	\bibliographystyle{IEEEtran}  
	\bibliography{MCCpaper}
}

\appendix
\subsection{Derivation of $\sigma_{j,k}^2$}

Let $\text{Si}(x) = \int_0^x \frac{\sin t}{t} dt$ and  $\text{Ci}(x) = \int_0^x \frac{\cos t}{t} dt$, then we have
\begin{equation}
	\begin{aligned}
		Q_{n, l} &
		= \frac{1}{(T_{l}^3) \pi^2} \int_{ (-B_{\sf e}/2-n\Delta_{l, \sf f})T_{l}}^{(B_{\sf e}/2-n\Delta_{l, \sf f})T_{l}} (q +  n\Delta_{l, \sf f}T_{l})^2 \text{sinc}^2(q) dq \\
		&= \frac{1}{(T_{l}^3 \pi^2)} \left( b_{n,l,0} + b_{n,l,1} + b_{n,l,2} \right), \nonumber
	\end{aligned}
\end{equation}
where $b_{n,l,0}, b_{n,l,1}$, and $b_{n,l,2}$ are given in \eqref{eqlong_7}.
\begin{figure*}
	\hrulefill
\begin{equation}
	\begin{aligned}
		b_{n,l,0} =& \int_{ (-B_{\sf e}/2-n\Delta_{l, \sf f})T_{l}}^{(B_{\sf e}/2-n\Delta_{l, \sf f})T_{l}}  \text{sin}^2(\pi q) dq =  B_{\sf e}T_l/2  - \frac{\sin (2\pi(B_{\sf e}/2-n\Delta_{l, \sf f})T_{l}) - \sin (2\pi(-B_{\sf e}/2-n\Delta_{l, \sf f})T_{l})}{4\pi}, \\
		b_{n,l,1} =&  \int_{ (-B_{\sf e}/2-n\Delta_{l, \sf f})T_{l}}^{(B_{\sf e}/2-n\Delta_{l, \sf f})T_{l}}  n\Delta_{l, \sf f}T_{l} \frac{(1- \cos(2\pi q))}{q} dq
		=n\Delta_{l, \sf f}T_{l} \left( \log\left[\left(B_{\sf e}/2-n\Delta_{l, \sf f} \right)T_{l} \right] - \log\left[\left(B_{\sf e}/2+n\Delta_{l, \sf f} \right)T_{l} \right]  \right)   \\
		&\hspace{5 cm}+ n\Delta_{l, \sf f}T_{l}\text{Ci}\left( 2 \pi (-B_{\sf e}/2-n\Delta_{l, \sf f})T_{l} \right) 
		- n\Delta_{l, \sf f}T_{l}\text{Ci}\left( 2\pi (B_{\sf e}/2-n\Delta_{l, \sf f})T_{l} \right), \\
		b_{n,l,2} =&\int_{ (-B_{\sf e}/2-n\Delta_{l, \sf f})T_{l}}^{(B_{\sf e}/2-n\Delta_{l, \sf f})T_{l}}  (n\Delta_{l, \sf f}T_{l})^2 \frac{\text{sin}^2(\pi q)}{q^2} dq
		=(n\Delta_{l, \sf f}T_{l})^2 \pi \left( \text{Si}\left( 2\pi(B_{\sf e}/2-n\Delta_{l, \sf f})T_{l} \right) + \text{Si}\left( 2\pi(B_{\sf e}/2+n\Delta_{l, \sf f})T_{l} \right) \right)\\
		&\hspace{5 cm}-(n\Delta_{l, \sf f}T_{l})^2\frac{\sin^2 (\pi(B_{\sf e}/2-n\Delta_{l, \sf f})T_{l} )}{(B_{\sf e}/2-n\Delta_{l, \sf f})T_{l} }-(n\Delta_{l, \sf f}T)^2\frac{\sin^2 (\pi(-B_{\sf e}/2-n\Delta_{l, \sf f})T_{l} )}{(B_{\sf e}/2+n\Delta_{l, \sf f})T_{l} }, \label{eqlong_7}
	\end{aligned}
\end{equation}
\hrulefill
\end{figure*}

On the other hand, $C_{j',i,l,n}^{j, i_{j},l_{j}, n_{j}}$ is rewritten as 
\begin{equation}
	\begin{aligned}
		C_{j',i,l,n}^{j, i_{j},l_{j}, n_{j}} 
		= \frac{\sum_{z=0}^4 c_{j',i,l,n,z}^{j, i_{j},l_{j}, n_{j}} }{{4\pi^4 T_{l}^2 T_{l_{j}}^2}}, \nonumber
	\end{aligned}
\end{equation}
where $c_{j',i,l,n,z}^{j, i_{j},l_{j}, n_{j}}, z=0,1,2,3,4$ are given in \eqref{eqlong_13a}-\eqref{eqlong_13e}, respectively. In order to obtain $c_{j',i,l,n,z}^{j, i_{j},l_{j}, n_{j}}, z=0,1,2,3,4$, we consider two different cases when $n_{j}\Delta_{l_{j}, \sf f}\neq  n\Delta_{l, \sf f}$ and $n_{j}\Delta_{l_{j}, \sf f} =  n\Delta_{l, \sf f}$.
\begin{figure*}
	\begin{subequations}
		\begin{alignat}{4}
			c_{j',i,l,n,0}^{j, i_{j},l_{j}, n_{j}} &=\int_{-B_{\sf e}/2}^{B_{\sf e}/2}  \frac{f^2   }{(f-n\Delta_{l, \sf f})^2(f-n_{j}\Delta_{l_{j}, \sf f})^2}df, \label{eqlong_13a} \\		
			c_{j',i,l,n,1}^{j, i_{j},l_{j}, n_{j}} &= \int_{-B_{\sf e}/2}^{B_{\sf e}/2}  \frac{f^2  \cos\left[2\pi(f(T_{l}- T_{l_{j}}) - n\Delta_{l, \sf f}T_{l} +n_{j}\Delta_{l_{j}, \sf f}T_{l_{j}} ) \right] }{2(f-n\Delta_{l, \sf f})^2(f-n_{j}\Delta_{l_{j}, \sf f})^2}df, \label{eqlong_13b} \\
			c_{j',i,l,n,2}^{j, i_{j},l_{j}, n_{j}} &=\int_{-B_{\sf e}/2}^{B_{\sf e}/2}  \frac{f^2  \cos\left[2\pi(f(T_{l}+ T_{l_{j}}) - n\Delta_{l, \sf f}T_{l} -n_{j}\Delta_{l_{j}, \sf f}T_{l_{j}} ) \right] }{2(f-n\Delta_{l, \sf f})^2(f-n_{j}\Delta_{l_{j}, \sf f})^2}df, 
			\label{eqlong_13c}\\
			c_{j',i,l,n,3}^{j, i_{j},l_{j}, n_{j}}&= -\int_{-B_{\sf e}/2}^{B_{\sf e}/2}  \frac{f^2  \cos\left[2\pi fT_{l} - 2\pi n \Delta_{l, \sf f}T_{l}  \right] }{(f-n\Delta_{l, \sf f})^2(f-n_{j}\Delta_{l_{j}, \sf f})^2}df, 
			\label{eqlong_13d}\\
			c_{j',i,l,n,4}^{j, i_{j},l_{j}, n_{j}}&=-\int_{-B_{\sf e}/2}^{B_{\sf e}/2}  \frac{f^2  \cos\left[2\pi fT_{l_{j}} - 2\pi n_{j}  \Delta_{l_{j}, \sf f}T_{l_{j}}  \right] }{(f-n\Delta_{l, \sf f})^2(f-n_{j}\Delta_{l_{j}, \sf f})^2}df. \label{eqlong_13e}
		\end{alignat}
	\end{subequations} 
\hrulefill
\end{figure*}

If $n_{j}\Delta_{l_{j}, \sf f}\neq  n\Delta_{l, \sf f}$, $c_{j',i,l,n,0}^{j, i_{j},l_{j}, n_{j}}$ is given in \eqref{eqlong_14},
\begin{figure*}
\begin{subequations}
	\begin{alignat}{4}
		\mathcal{G}_1(n\Delta_{l, \sf f}, n_{j}\Delta_{l_{j}, \sf f}) =& \frac{d_{n,l,0}^{n_j, l_j}   }{ (f-n\Delta_{l, \sf f})} +   \frac{d_{n,l,1}^{n_j, l_j}   }{ (f-n\Delta_{l, \sf f})^2} +\frac{d_{n,l,2}^{n_j, l_j}   }{ (f-n_{j}\Delta_{l_{j}, \sf f})}+ \frac{d_{n,l,3}^{n_j, l_j}   }{ (f-n_{j}\Delta_{l_{j}, \sf f})^2},\\
		\mathcal{G}_2(n\Delta_{l, \sf f}) =&  \frac{d_{n,l,4}^{n_j, l_j}   }{ (f-n\Delta_{l, \sf f})^2} +   \frac{d_{n,l,5}^{n_j, l_j}   }{ (f-n\Delta_{l, \sf f})^3} +\frac{d_{n,l,6}^{n_j, l_j}   }{ (f-n\Delta_{l, \sf f})^4},\\
		c_{j',i,l,n,0}^{j, i_{j},l_{j}, n_{j}} 
		=& \int_{-B_{\sf e}/2}^{B_{\sf e}/2} \mathcal{G}_1(n\Delta_{l, \sf f}, n_{j}\Delta_{l_{j}, \sf f}) df \nonumber\\
		=& d_{n,l,0}^{n_j, l_j}\ln (|B_{\sf e}-n\Delta_{l, \sf f}|) - d_{n,l,0}^{n_j, l_j}\ln (|B_{\sf e}+n\Delta_{l, \sf f}|)-\frac{d_{n,l,1}^{n_j, l_j}}{B_{\sf e}-n\Delta_{l, \sf f}} - \frac{d_{n,l,1}^{n_j, l_j}}{B_{\sf e}+n\Delta_{l, \sf f}} \nonumber\\
		&+ d_{n,l,2}^{n_j, l_j}\ln (|B_{\sf e}-n_{j}\Delta_{l_{j}, \sf f}|) - d_{n,l,2}^{n_j, l_j}\ln (|B_{\sf e}+n_{j}\Delta_{l_{j}, \sf f}|)-\frac{d_{n,l,3}^{n_j, l_j}}{B_{\sf e}-n_{j}\Delta_{l_{j}, \sf f}} - \frac{d_{n,l,3}^{n_j, l_j}}{B_{\sf e}+n_{j}\Delta_{l_{j}, \sf f}}, \text{ if } n_{j}\Delta_{l_{j}, \sf f}\neq  n\Delta_{l, \sf f}, \label{eqlong_14} \\
			c_{j',i,l,n,0}^{j, i_{j},l_{j}, n_{j}} 
			=& \int_{-B_{\sf e}/2}^{B_{\sf e}/2} \mathcal{G}_2(n\Delta_{l, \sf f}) df  =\frac{d_{n,l,4}^{n_j, l_j}}{B_{\sf e}-n\Delta_{l, \sf f}} - \frac{d_{n,l,4}^{n_j, l_j}}{B_{\sf e}+n\Delta_{l, \sf f}} - \frac{d_{n,l,5}^{n_j, l_j}}{2(B_{\sf e}-n\Delta_{l, \sf f})^2} - \frac{d_{n,l,5}^{n_j, l_j}}{2(B_{\sf e}+n\Delta_{l, \sf f})^2} \nonumber\\
			 & \hspace{3.3 cm} -\frac{d_{n,l,6}^{n_j, l_j}}{2(B_{\sf e}-n\Delta_{l, \sf f})^3} - \frac{d_{n,l,6}^{n_j, l_j}}{3(B_{\sf e}+n\Delta_{l, \sf f})^3}, \text{ if } n_{j}\Delta_{l_{j}, \sf f}\ =  n\Delta_{l, \sf f}. \label{eqlong_15}
		\end{alignat}
	\end{subequations} 
\end{figure*}
where $d_{n,l,0}^{n_j, l_j}, d_{n,l,1}^{n_j, l_j}, d_{n,l,2}^{n_j, l_j}, d_{n,l,3}^{n_j, l_j}$ are the functions of $n\Delta_{l, \sf f}$ and $n_{j}\Delta_{l_{j}, \sf f}$ which is the root of the equations in \eqref{eq45}.
\begin{equation}
	\begin{cases}
		&d_{n,l,0}^{n_j, l_j} + d_{n,l,2}^{n_j, l_j} = 0,\\
		&d_{n,l,0}^{n_j, l_j}(-n_{j}\Delta_{l_{j}, \sf f}+n\Delta_{l, \sf f}) + d_{n,l,1}^{n_j, l_j} + d_{n,l,3}^{n_j, l_j} = 1, \\
		&d_{n,l,0}^{n_j, l_j}(n_{j}^2\Delta_{l_{j}, \sf f}^2 - n^2\Delta_{l, \sf f}^2) - 2 d_{n,l,1}^{n_j, l_j} n_{j}\Delta_{l_{j}, \sf f} \\
		& \hspace{ 4.9 cm} - 2d_{n,l,3}^{n_j, l_j}n\Delta_{l, \sf f} = 0,\\
		& d_{n,l,0}^{n_j, l_j}(-n\Delta_{l, \sf f}n_{j}^2\Delta_{l_{j}, \sf f}^2 +n^2\Delta_{l, \sf f}^2n_{j}\Delta_{l_{j}, \sf f} ) \\
		&\hspace{ 2.6 cm} + d_{n,l,1}^{n_j, l_j} n_{j}^2\Delta_{l_{j}, \sf f}^2 + d_{n,l,3}^{n_j, l_j} n^2\Delta_{l, \sf f}^2 = 0. \label{eq45}
	\end{cases}
\end{equation}

If $n_{j}\Delta_{l_{j}, \sf f}= n\Delta_{l, \sf f}$,  $c_{j',i,l,n,0}^{j, i_{j},l_{j}, n_{j}}$ is given in \eqref{eqlong_15},
where $d_{n,l,4}^{n_j, l_j} = 1, d_{n,l,5}^{n_j, l_j} = 2n\Delta_{l, \sf f}$, and $d_{n,l,6}^{n_j, l_j}=n^2\Delta_{l, \sf f}^2$.

\begin{figure*}
	\begin{equation}
		\begin{aligned}
			c_{j',i,l,n,1}^{j, i_{j},l_{j}, n_{j}} 
			=&\frac{1}{2} \int_{-B_{\sf e}/2}^{B_{\sf e}/2} \mathcal{G}_1(n\Delta_{l, \sf f}, n_{j}\Delta_{l_{j}, \sf f}) \cos\left[2\pi(f(T_{l}- T_{l_{j}}) - n\Delta_{l, \sf f}T_{l} +n_{j}\Delta_{l_{j}, \sf f}T_{l_{j}} ) \right]df\\
			=&  \frac{d_{n,l,0}^{n_j, l_j}}{2} \mathcal{E}_1(B_{\sf e},2\pi f(T_{l}- T_{l_{j}}), - 2\pi n\Delta_{l, \sf f}T_{l} + 2\pi n_{j}\Delta_{l_{j}, \sf f}T_{l_{j}}, n\Delta_{l, \sf f}) \\
			&+ \frac{d_{n,l,1}^{n_j, l_j}}{2} \mathcal{E}_2(B_{\sf e},2\pi f(T_{l}- T_{l_{j}}), - 2\pi n\Delta_{l, \sf f}T_{l} + 2\pi n_{j}\Delta_{l_{j}, \sf f}T_{l_{j}}, n\Delta_{l, \sf f})\\
			&+ \frac{d_{n,l,2}^{n_j, l_j}}{2} \mathcal{E}_1(B_{\sf e},2\pi f(T_{l}- T_{l_{j}}), - 2\pi n\Delta_{l, \sf f}T_{l} + 2\pi n_{j}\Delta_{l_{j}, \sf f}T_{l_{j}}, n_{j}\Delta_{l_{j}, \sf f})\\
			&+ \frac{d_{n,l,3}^{n_j, l_j}}{2} \mathcal{E}_2(B_{\sf e},2\pi f(T_{l}- T_{l_{j}}), - 2\pi n\Delta_{l, \sf f}T_{l} + 2\pi n_{j}\Delta_{l_{j}, \sf f}T_{l_{j}}, n_{j}\Delta_{l_{j}, \sf f}), \text{ if } n_{j}\Delta_{l_{j}, \sf f}\neq  n\Delta_{l, \sf f},\\
			c_{j',i,l,n,1}^{j, i_{j},l_{j}, n_{j}} 
			=&\frac{1}{2} \int_{-B_{\sf e}/2}^{B_{\sf e}/2} \mathcal{G}_2(n\Delta_{l, \sf f}, n_{j}\Delta_{l_{j}, \sf f}) \cos\left[2\pi(f(T_{l}- T_{l_{j}}) - n\Delta_{l, \sf f}T_{l} +n_{j}\Delta_{l_{j}, \sf f}T_{l_{j}} ) \right]df\\
			=&  \frac{d_{n,l,4}^{n_j, l_j}}{2} \mathcal{E}_2(B_{\sf e},2\pi f(T_{l}- T_{l_{j}}), - 2\pi n\Delta_{l, \sf f}T_{l} + 2\pi n\Delta_{l, \sf f}T_{l_{j}}, n\Delta_{l, \sf f}) \\
			&+ \frac{d_{n,l,5}^{n_j, l_j}}{2} \mathcal{E}_3(B_{\sf e},2\pi f(T_{l}- T_{l_{j}}), - 2\pi n\Delta_{l, \sf f}T_{l} + 2\pi n\Delta_{l, \sf f}T_{l_{j}}, n\Delta_{l, \sf f})\\
			&+ \frac{d_{n,l,6}^{n_j, l_j}}{2} \mathcal{E}_4(B_{\sf e},2\pi f(T_{l}- T_{l_{j}}), - 2\pi n\Delta_{l, \sf f}T_{l} + 2\pi n\Delta_{l, \sf f}T_{l_{j}}, n\Delta_{l, \sf f}), \text{ if } n_{j}\Delta_{l_{j}, \sf f}= n\Delta_{l, \sf f}. \label{eqlong_21}
		\end{aligned}
	\end{equation}
\end{figure*}

\begin{figure*}
	\begin{equation}
		\begin{aligned}
			c_{j',i,l,n,2}^{j, i_{j},l_{j}, n_{j}} 
			=&\frac{1}{2} \int_{-B_{\sf e}/2}^{B_{\sf e}/2} \mathcal{G}_1(n\Delta_{l, \sf f}, n_{j}\Delta_{l_{j}, \sf f}) \cos\left[2\pi(f(T_{l}+T_{l_{j}}) - n\Delta_{l, \sf f}T_{l} +n_{j}\Delta_{l_{j}, \sf f}T_{l_{j}} ) \right]df \\
			=&  \frac{d_{n,l,0}^{n_j, l_j}}{2} \mathcal{E}_1(B_{\sf e},2\pi f(T_{l}+T_{l_{j}}), - 2\pi n\Delta_{l, \sf f}T_{l} - 2\pi n_{j}\Delta_{l_{j}, \sf f}T_{l_{j}}, n\Delta_{l, \sf f}) \\
			&+ \frac{d_{n,l,1}^{n_j, l_j}}{2} \mathcal{E}_2(B_{\sf e},2\pi f(T_{l}+T_{l_{j}}), - 2\pi n\Delta_{l, \sf f}T_{l} - 2\pi n_{j}\Delta_{l_{j}, \sf f}T_{l_{j}}, n\Delta_{l, \sf f})\\
			&+ \frac{d_{n,l,2}^{n_j, l_j}}{2} \mathcal{E}_1(B_{\sf e},2\pi f(T_{l}+T_{l_{j}}), - 2\pi n\Delta_{l, \sf f}T_{l} - 2\pi n_{j}\Delta_{l_{j}, \sf f}T_{l_{j}}, n_{j}\Delta_{l_{j}, \sf f})\\
			&+ \frac{d_{n,l,3}^{n_j, l_j}}{2} \mathcal{E}_2(B_{\sf e},2\pi f(T_{l}+T_{l_{j}}), - 2\pi n\Delta_{l, \sf f}T_{l} - 2\pi n_{j}\Delta_{l_{j}, \sf f}T_{l_{j}}, n_{j}\Delta_{l_{j}, \sf f}), \text{ if } n_{j}\Delta_{l_{j}, \sf f}\neq  n\Delta_{l, \sf f},\\
			c_{j',i,l,n,2}^{j, i_{j},l_{j}, n_{j}} 
			=&\frac{1}{2} \int_{-B_{\sf e}/2}^{B_{\sf e}/2} \mathcal{G}_2(n\Delta_{l, \sf f}, n_{j}\Delta_{l_{j}, \sf f}) \cos\left[2\pi(f(T_{l}+T_{l_{j}}) - n\Delta_{l, \sf f}T_{l} +n_{j}\Delta_{l_{j}, \sf f}T_{l_{j}} ) \right]df \\
			=&  \frac{d_{n,l,4}^{n_j, l_j}}{2} \mathcal{E}_2(B_{\sf e},2\pi f(T_{l}+T_{l_{j}}), - 2\pi n\Delta_{l, \sf f}T_{l} + 2\pi n\Delta_{l, \sf f}T_{l_{j}}, n\Delta_{l, \sf f}) \\
			&+ \frac{d_{n,l,5}^{n_j, l_j}}{2} \mathcal{E}_3(B_{\sf e},2\pi f(T_{l}+T_{l_{j}}), - 2\pi n\Delta_{l, \sf f}T_{l} + 2\pi n\Delta_{l, \sf f}T_{l_{j}}, n\Delta_{l, \sf f})\\
			&+ \frac{d_{n,l,6}^{n_j, l_j}}{2} \mathcal{E}_4(B_{\sf e},2\pi f(T_{l}+T_{l_{j}}), - 2\pi n\Delta_{l, \sf f}T_{l} + 2\pi n\Delta_{l, \sf f}T_{l_{j}}, n\Delta_{l, \sf f}), \text{ if } n_{j}\Delta_{l_{j}, \sf f}= n\Delta_{l, \sf f}. \label{eqlong_22}
		\end{aligned}
	\end{equation}
\end{figure*}

\begin{figure*}
	\begin{equation}
		\begin{aligned}
			c_{j',i,l,n,3}^{j, i_{j},l_{j}, n_{j}} 
			=& \int_{-B_{\sf e}/2}^{B_{\sf e}/2} \mathcal{G}_1(n\Delta_{l, \sf f}, n_{j}\Delta_{l_{j}, \sf f}) \cos\left[2\pi fT_{l} - 2\pi n\Delta_{l, \sf f}T_{l}  \right]df \\
			=&  {d_{n,l,0}^{n_j, l_j}} \mathcal{E}_1(B_{\sf e},2\pi fT_{l}, - 2\pi n\Delta_{l, \sf f}T_{l}, n\Delta_{l, \sf f}) + {d_{n,l,1}^{n_j, l_j}} \mathcal{E}_2(B_{\sf e},2\pi fT_{l}, - 2\pi n\Delta_{l, \sf f}T_{l}, n\Delta_{l, \sf f})\\
			&+ {d_{n,l,2}^{n_j, l_j}} \mathcal{E}_1(B_{\sf e},2\pi fT_{l}, - 2\pi n\Delta_{l, \sf f}T_{l}, n_{j}\Delta_{l_{j}, \sf f})+ {d_{n,l,3}^{n_j, l_j}} \mathcal{E}_2(B_{\sf e},2\pi fT_{l}, - 2\pi n\Delta_{l, \sf f}T_{l}, n_{j}\Delta_{l_{j}, \sf f}), \text{ if } n_{j}\Delta_{l_{j}, \sf f}\neq  n\Delta_{l, \sf f},\\
			c_{j',i,l,n,3}^{j, i_{j},l_{j}, n_{j}} 
			=& \int_{-B_{\sf e}/2}^{B_{\sf e}/2} \mathcal{G}_2(n\Delta_{l, \sf f}, n_{j}\Delta_{l_{j}, \sf f}) \cos\left[2\pi fT_{l} - 2\pi n\Delta_{l, \sf f}T_{l}  \right]df \\
			=&  {d_{n,l,4}^{n_j, l_j}} \mathcal{E}_2(B_{\sf e},2\pi fT_{l}, - 2\pi n\Delta_{l, \sf f}T_{l}, n\Delta_{l, \sf f}) + {d_{n,l,5}^{n_j, l_j}} \mathcal{E}_3(B_{\sf e},2\pi fT_{l}, - 2\pi n\Delta_{l, \sf f}T_{l}, n\Delta_{l, \sf f})\\
			&+ {d_{n,l,6}^{n_j, l_j}} \mathcal{E}_4(B_{\sf e},2\pi fT_{l}, - 2\pi n\Delta_{l, \sf f}T_{l}, n\Delta_{l, \sf f}), \text{ if } n_{j}\Delta_{l_{j}, \sf f}= n\Delta_{l, \sf f}. \label{eqlong_23}
		\end{aligned}
	\end{equation}
\end{figure*}

\begin{figure*}
	\begin{equation}
		\begin{aligned}
			c_{j',i,l,n,4}^{j, i_{j},l_{j}, n_{j}} 
			=& \int_{-B_{\sf e}/2}^{B_{\sf e}/2} \mathcal{G}_1(n_j\Delta_{l_j, \sf f}, n_{j}\Delta_{l_{j}, \sf f}) \cos\left[2\pi fT_{l_j} - 2\pi n_j\Delta_{l_j, \sf f}T_{l_j} \right]df \\
			=&  {d_{n,l,0}^{n_j, l_j}} \mathcal{E}_1(B_{\sf e},2\pi fT_{l_j}, - 2\pi n_j\Delta_{l_j, \sf f}T_{l_j}, n_j\Delta_{l_j, \sf f}) + {d_{n,l,1}^{n_j, l_j}} \mathcal{E}_2(B_{\sf e},2\pi fT_{l_j}, - 2\pi n_j\Delta_{l_j, \sf f}T_{l_j}, n_j\Delta_{l_j, \sf f})\\
			&+ {d_{n,l,2}^{n_j, l_j}} \mathcal{E}_1(B_{\sf e},2\pi fT_{l_j}, - 2\pi n_j\Delta_{l_j, \sf f}T_{l_j}, n_{j}\Delta_{l_{j}, \sf f}) \\
			&+ {d_{n,l,3}^{n_j, l_j}} \mathcal{E}_2(B_{\sf e},2\pi fT_{l_j}, - 2\pi n_j\Delta_{l_j, \sf f}T_{l_j}, n_{j}\Delta_{l_{j}, \sf f}), \text{ if } n_{j}\Delta_{l_{j}, \sf f}\neq  n_j\Delta_{l_j, \sf f},\\
			c_{j',i,l,n,4}^{j, i_{j},l_{j}, n_{j}} 
			=& \int_{-B_{\sf e}/2}^{B_{\sf e}/2} \mathcal{G}_2(n_j\Delta_{l_j, \sf f}, n_{j}\Delta_{l_{j}, \sf f}) \cos\left[2\pi fT_{l_j} - 2\pi n_j\Delta_{l_j, \sf f}T_{l_j} \right]df \\
			=&  {d_{n,l,4}^{n_j, l_j}} \mathcal{E}_2(B_{\sf e},2\pi fT_{l_j}, - 2\pi n_j\Delta_{l_j, \sf f}T_{l_j}, n_j\Delta_{l_j, \sf f}) + {d_{n,l,5}^{n_j, l_j}} \mathcal{E}_3(B_{\sf e},2\pi fT_{l_j}, - 2\pi n_j\Delta_{l_j, \sf f}T_{l_j}, n_j\Delta_{l_j, \sf f})\\
			&+ {d_{n,l,6}^{n_j, l_j}} \mathcal{E}_4(B_{\sf e},2\pi fT_{l_j}, - 2\pi n_j\Delta_{l_j, \sf f}T_{l_j}, n_j\Delta_{l_j, \sf f}), \text{ if } n_{j}\Delta_{l_{j}, \sf f}= n_j\Delta_{l_j, \sf f}. \label{eqlong_24}
		\end{aligned}
	\end{equation}
\end{figure*}

	Similarly, $c_{j',i,l,n,z}^{j, i_{j},l_{j}, n_{j}}, z=1,2,3,4$ are computed as in \eqref{eqlong_21}, \eqref{eqlong_21}, \eqref{eqlong_21}, \eqref{eqlong_21}, respectively, where $\mathcal{E}_i(B_{\sf e},a,b,c), i=1,2,3,4$ are given in \eqref{eqlong_30}, \eqref{eqlong_31}, \eqref{eqlong_32}, \eqref{eqlong_33}, \eqref{eqlong_34}, respectively.

\begin{figure*}
	\begin{subequations}
		\begin{alignat}{4}
			\mathcal{E}_0(B_{\sf e},a,b,c) =& \int_{-B_{\sf e}/2}^{B_{\sf e}/2} \frac{\sin(ax +b)}{(x-c)}dx =a \int_{-B_{\sf e}/2 - c}^{B_{\sf e}/2-c} \frac{\sin(ax) \cos(ac +b) + \cos(ax) \sin(ac +b)}{ax}dx \nonumber \\
			=& a\cos(ac +b) \text{Si}(ax)|_{-B_{\sf e}/2 - c}^{B_{\sf e}/2-c} + a\sin(ac +b) \text{Ci}(ax)|_{-B_{\sf e}/2 - c}^{B_{\sf e}/2-c}, \label{eqlong_30}\\
			\mathcal{E}_1(B_{\sf e},a,b,c) =& \int_{-B_{\sf e}/2}^{B_{\sf e}/2} \frac{\cos(ax +b)}{(x-c)}dx =
			a \int_{-B_{\sf e}/2 - c}^{B_{\sf e}/2-c} \frac{\cos(ax) \cos(ac +b) - \sin(ax) \sin(ac +b)}{ax}dx \nonumber \\
			=& a\cos(ac +b) \text{Ci}(ax)|_{-B_{\sf e}/2 - c}^{B_{\sf e}/2-c} - a\sin(ac +b) \text{Si}(ax)|_{-B_{\sf e}/2 - c}^{B_{\sf e}/2-c}, \label{eqlong_31}\\
			\mathcal{E}_2(B_{\sf e},a,b,c) =& \int_{-B_{\sf e}/2}^{B_{\sf e}/2} \frac{\cos(ax +b)}{(x-c)^2}dx = -\frac{\cos(ax +b)}{(x-c)}|_{-B_{\sf e}/2}^{B_{\sf e}/2} - a\mathcal{E}_0(B_{\sf e},a,b,c), \label{eqlong_32}\\
			\mathcal{E}_3(B_{\sf e},a,b,c)		=& \int_{-B_{\sf e}/2}^{B_{\sf e}/2} \frac{\cos(ax +b)}{(x-c)^3}dx
			=-\frac{\cos(ax+b)}{2(x-c)^2} \Big|_{-B_{\sf e}/2}^{B_{\sf e}/2} + \frac{a\sin(ax+b)}{2(x-c)}\Big|_{-B_{\sf e}/2}^{B_{\sf e}/2} -\frac{a^2}{2} \mathcal{E}_1(B_{\sf e},a,b,c), \label{eqlong_33}\\
			\mathcal{E}_4(B_{\sf e},a,b,c)	=& \int_{-B_{\sf e}/2}^{B_{\sf e}/2} \frac{\cos(ax +b)}{(x-c)^4}dx \nonumber \\
			=&-\frac{\cos(ax+b)}{3(x-c)^3} \Big|_{-B_{\sf e}/2}^{B_{\sf e}/2} + \frac{a\sin(ax+b)}{6(x-c)^2}\Big|_{-B_{\sf e}/2}^{B_{\sf e}/2}  +\frac{a^2\cos(ax+b)}{6(x-c)} \Big|_{-B_{\sf e}/2}^{B_{\sf e}/2} +\frac{a^3}{6} \mathcal{E}_0(B_{\sf e},a,b,c). \label{eqlong_34}
		\end{alignat}
	\end{subequations} 
	\hrulefill
\end{figure*}

\end{document}